\documentclass{article}
\usepackage{amssymb,amsmath,mathrsfs}
\usepackage[colorlinks=true, pdfstartview=FitV, linkcolor=blue, citecolor=blue, urlcolor=blue]{hyperref}

\usepackage{geometry} 
\usepackage{framed,color}
\definecolor{shadecolor}{rgb}{0.9, 0.9, 0.81}
\usepackage{ tikz, pgfplots, todonotes,amssymb}
\usetikzlibrary{decorations.markings,arrows, calc, shapes, intersections, patterns}

\def \scr{\mathscr}

\def \div{ {\rm div}\,}

\def \DD{\mathbb D}

\def \CC{{\mathcal C}}
\def \ds{\displaystyle}
\def\le{\left}
\def \bs{\boldsymbol}
\def \QED{\hfill $\blacksquare$\par \vskip 4pt}
\def\ri{\right}

\def\bea#1\eea{\begin{align}#1\end{align}}
\def \wt{ \widetilde }
\def \&{\hspace{-15pt}&}
\def \d{{\mathrm d}}
\def\res{\mathop{\mathrm {res}}}
\def \wh{\widehat}

\def \bd#1\ed{ \begin{definition} #1 \hfill $\blacktriangle$\end{definition}}
\def \br#1\er {\begin{remark} {#1 \hfill $\blacktriangle$}\end{remark}}
\newtheorem{theorem}{Theorem}[section]
\newtheorem{definition}[theorem]{Definition}
\newtheorem{problem}[theorem]{Problem}

\newtheorem{proposition}[theorem]{Proposition}
\newtheorem{corollary}[theorem]{Corollary}
\newtheorem{remark}[theorem]{Remark}
\newtheorem{lemma}[theorem]{Lemma}
\newtheorem{example}[theorem]{Example}

\def \bp#1\ep{
 \definecolor{shadecolor}{rgb}{0.95, 0.95, 0.86}
 \begin{shaded}\begin{proposition} #1 \end{proposition}\end{shaded}}
\def \bc#1\ec{
 \definecolor{shadecolor}{rgb}{0.95, 0.95, 0.86}
 \begin{shaded}\begin{corollary}#1\end{corollary}\end{shaded}}

\def \bt#1\et{
 \definecolor{shadecolor}{rgb}{0.95, 0.95, 0.86}
 \begin{shaded}\begin{theorem} #1 \end{theorem}\end{shaded}}

\def \bl#1\el{
 \definecolor{shadecolor}{rgb}{0.95, 0.95, 0.86}
 \begin{shaded}\begin{lemma} #1 \end{lemma}\end{shaded}}

\def \L {\mathcal L}
\def\K{\mathcal K}

\def\nn{\nonumber}
\def\be{\begin{equation}}
\def\ee{\end{equation}}
\def\I{\mathcal I}

\def\ba{\begin{array}}
\def\ea{\end{array}}
\makeatletter
\@addtoreset{equation}{section}
\makeatother

\def \eqref #1{(\ref{#1})}
\def \1{\mathbf 1}
\def\M{{\mathcal M}}
\def\C{{\mathbb C}}
\def\Z{{\mathbb  Z}}
\def\R{{\mathbb R}}

\def\l{\lambda}

\def\p{{\mathbf p}}

\def\N{{\mathbb N}}

\def\pa{\partial}
\def\ov{\overline}
\def \Re{\mathrm {Re}\,}
\def \Im {\mathrm {Im}\,}

\def \z{\zeta}

\begin{document}
\baselineskip 15pt plus 1pt minus 1pt

\vspace{0.2cm}
\begin{center}
\begin{Large}
\fontfamily{cmss}
\fontsize{17pt}{27pt}
\selectfont
\textbf{
The generalized Chebotarev problem in higher genus}
\end{Large}\\
\bigskip
\begin{large} {M.
Bertola}$^{\ddagger,\sharp}$\footnote{Marco.Bertola@concordia.ca}
\end{large}
\\
\bigskip
\begin{small}
$^{\ddagger}$ {\em Department of Mathematics and
Statistics, Concordia University\\ 1455 de Maisonneuve W., Montr\'eal, Qu\'ebec,
Canada H3G 1M8} \\
$^{\sharp}$ {\em Centre de recherches math\'ematiques, Universit\'e\ de
Montr\'eal } \\
\end{small}
\bigskip
{\bf Abstract}\end{center}
We consider the extension to higher genus Riemann surfaces of the classical Chebotarev problem, with a view towards the development of the theory of Pad\'e\ approximants on algebraic curves. 
To this end we define an appropriate notion of capacity that mimics the standard one, following works of Chirka and of the author and collaborators. 
The nontrivial topology of the Riemann surface requires further specification of the ``homotopy class'' of the continua in the solution of the Chebotarev problem. We also discuss the relationship of this problem to the theory of Jenkins-Strebel quadratic differentials.

\vspace{0.7cm}

\vskip 15pt
{}

\tableofcontents
\section{Introduction}
The theory of Pad\'e\ approximation is  a useful and beautiful branch of approximation theory, which touches upon several distant areas such as number theory, random matrices, integrable systems, to mention a few. 
By its nature the classical theory is formulated as an approximation tool for functions on the complex plane (or subsets thereof), namely on the simplest of Riemann surfaces, the Riemann sphere $\mathbb P^1$. 

In recent years we have started a  program aimed at generalizing  the theory and potentially its applications  to the case where the Riemann sphere is replaced by a more general Riemann surface $\CC$  of arbitrary genus. 
In \cite{BertoPade} we introduced the notion of Pad\'e approximants for an appropriate generalization of Stieltjes-Markov functions of measures on paths embedded in $\CC$, showing how the denominators, meromorphic functions with prescribed poles, satisfy an orthogonality that mimics the relationship between orthogonal polynomials and measures in the complex plane.  In particular it was shown that the  denominators and suitable Cauchy transforms thereof can be arranged in a $2\times 2$ matrix  that satisfies a Riemann-Hilber problem, much like ordinary orthogonal polynomials do \cite{FIK}. This allowed  \cite{BertoEllOPs} to analyze their behaviour for large degree using the nonlinear steepest descent method much like  \cite{FIK, DKMVZ}. Another application to matrix orthogonal polynomials was developed in \cite{BertoAbelianization}. 

This paper is then born out of a necessity that arises in the theory of nonlinear steepest descent method, or Deift-Zhou method \cite{DZ}, when   generalized to Riemann-Hilbert problems on higher genus Riemann surfaces. Namely, that of constructing an appropriate equilibrium measure, which is then used to start the Deift-Zhou method for analyzing the large degree asymptotics of orthogonal rational functions. 

It is  well known \cite{Stahl, Stahl2,Stahl3} that the zeros of denominators  of the diagonal  Pad\'e approximants  for a germ of  analytic function $F(z)$ near $z=\infty$  accumulate on a set $\frak F$ with the property that 
\begin{itemize}
\item  the function $F$ admits analytic extension to $\mathbb P^1\setminus \frak F$;
\item the set $\frak F$ is the unique set of minimal logarithmic capacity amongst all sets satisfying the previous condition. 
\end{itemize}
The typical situation is that of a function $F$ admitting unrestricted analytic continuation to the universal cover of $\mathbb P^1$ minus a finite set, $E$, of branch-points. The reader can think of an algebraic function $F$ for example.  In this context the function $F$ does not play a particularly prominent role; just the set of branchpoints $E$ and the overall topology of the maximal set of analytic continuation do. 

If, for example, $F$ can only be extended to a simply connected domain containing $\infty$, then $\frak F$ must be a continuum (connected compact  set) containing the branch-points $E$. Then this can be reformulated as the classical Chebotarev problem. 
\begin {problem}[P\'olya in a letter to Chebotarev \cite{Polya}, solved by Grotzsch \cite{Grotzsch} and Lavrentieff \cite{Lav2}]
For a given finite set $E=\{e_1,\dots, e_N\}\subset \C$, find a continuum $\frak F$ containing $E$ and of minimal logarithmic capacity amongst all such sets.
\end{problem}
We will call $E$ the {\it anchor set} and $e\in E$ an {\it anchor}. 
It may well be that $F$ can be extended analytically to a set that is not necessarily simply connected. The ``connectivity'' of such a set then plays an important role. For example, the function 
$$
F(z) = \sqrt[4]{\frac {z-2}{z-1}} + \sqrt[3]{\frac {z+i}{z-i}}
$$
clearly  has four branchpoints $E = \{1,2,i,-i\} $ but can be extended to any set where the branchcuts extends from $1$ to $2$ and from $i$ to $-i$: in other words it is not necessary for the domain of analyticity to be simply connected. Thus the set of minimal capacity $\frak F$ in Stahl's theory needs not be a continuum containing the anchors.  It must however be a {\it poly-continuum} (finite union of continua)  where $1,2$ belong to the same connected component and so do $i,-i$, but not necessarily the same component for the two pairs.  This leads to a {\it generalized} Chebotarev problem, where the class over which we minimize the capacity is that of poly-continua with the required ``connectivity pattern'', see also \cite{MFRakh}.  A numerical approach to generalized Chebotarev poly-continua in the plane can be found in \cite{BertoSAPM}. 
\paragraph{Summary of results.}
We can easily imagine the generalization of this problem to an arbitrary Riemann surface $\CC$ of genus $g$, and indeed this is the focus of this paper. 

Suppose now $E\subset \CC$ is a finite set of anchors on $\CC$ and we are seeking, similarly to the previously described situation, a  poly-continuum $\K$ of minimal capacity that preserves an appropriate  ``connectivity pattern'', $\scr P$:
the precise meaning is given in Section \ref{connpatt}, see Definition \ref{defclass}. Briefly,  it is a prescription of the homotopy classes (in the relative homotopy $\pi_1(\CC,E)$) that  the poly-continuum defines by connecting the anchors together. 

The definition of capacity is relatively simple and we follow ideas in  \cite{Chirka} and \cite{BertoGrootKuijlaars}. Like the ordinary logarithmic capacity in the complex plane, the notion requires to specify a reference point relative to which the capacity is computed. In the standard case this is usually the point at infinity and in our setting it is just a point on $\CC$ which we denote by $\infty$ as a reference to the standard case. 

Let us fix  a  point $\infty \in \CC$ and consider a compact set $\K\subset \CC\setminus \{\infty\}$  which is regular for the Dirichlet problem. Let $G_\K(p;\infty)$ denote the standard Green's function of the complement of $\K$ with logarithmic pole at the point $\infty$. Let us choose  a local coordinate $z_\infty$ near $\infty$ and   define  the capacity as 
\be
\ln {\rm Cap}_{z_\infty}(\K):= \lim_{p\to \infty} \le(\ln \frac 1{|z_\infty(p)|} -G_\K(p;\infty) \ri). 
\ee
 It  has the property that if $w_\infty$ is a different local coordinate then 
\be
\label{capzw}
{\rm Cap}_{z_\infty}(\K) = {\rm Cap}_{w_\infty}(\K)  \le|\frac {\d w_\infty}{\d z_\infty}\ri|_{p=\infty}.
\ee
The classical Frostman theorem relates the capacity of a set to a probability measure supported on its boundary, and Chirka \cite{Chirka} showed how to extend this result to Riemann surfaces. His very general results can be specialized also to the above notion of capacity. See Proposition \ref{propFrostman}.

In Section \ref{secChebo} we get to the first result. We first review Chirka's results that prove that the capacity ${\rm Cap}_{z_\infty}$ is continuous with respect to the Hausdorff topology on the class of poly-continua we are interested in, summarizing Chirka's result in  Lemma \ref{lemmacont}. This continuity allows us to prove Theorem \ref{thmP} which we rephrase below.

\begin{shaded}
{\bf Theorem}  For a given set of anchors $E$ and any connectivity pattern, there exists a unique poly-continuum that realizes the minimum of the capacity within this class. 
\end{shaded}
\noindent See Theorem \ref{thmP} for more precise formulation.  The existence is a simple consequence of the continuity of the capacity. The uniqueness is a slightly more delicate matter which leads us to the next section. 

In Section \ref{secminimal} we show how sets of minimal capacity are related to trajectories of quadratic differentials. The core ideas here are not new and have been used in \cite{Nuttall, GoncharRakhmanov, MFRakh}. One of the main features which we can state here in simple terms is the $S$-property introduced by Gonchar-Rakhmanov and Stahl \cite{GoncharRakhmanov, Stahl,Stahl2,Stahl3}. 
In short, the minimal capacity sets $\K$ are finite unions of analytic arcs, possibly meeting transversally at the endpoints, and the Green's function $G_\K(p)= G_\K(p;\infty)$ has the $S$-property
\be
\frac {\pa G_\K(p)}{\pa {\bf n}_+} = \frac {\pa G_\K(p)}{\pa {\bf n}_-}
\label{Spropertyintro},
\ee
where ${\bf n}_\pm$ denote the two opposite normals at a smooth point of $\K$. 
 Furthermore, using the tool of Schiffer variations \cite{SS54}  as in  \cite{RakhPer, MFRakh, KuijlaarsSilva, BertoGrootKuijlaars}we show that the expression $\mathfrak Q=(2\pa_p G_\K(p;\infty))^2$ (with $\pa_p$ the Wirthinger operator \eqref{Wirtinger}) extends to a meromorphic quadratic differential on $\CC$ such that: 
 \begin{itemize}
 \item $\frak Q$ has a double pole at $\infty$ with bi-residue $1$: the bi-residue of a quadratic differential at a double pole is the coefficient of the term $w^{-2}$ in any local coordinate. Thus  in any local coordinate $w$ near $\infty$ the differential has the expression
 \be
 \mathfrak Q(p) = \frac {\d w^2}{w^2} (1 + \mathcal O(w)).
 \label{biresidue}
 \ee
 \item $\frak Q$ has at most simple poles at the anchor set $E$.
 \item all contour integrals of $\sqrt{\frak Q}$ on the Riemann surface $\wh \CC$ defined by $v^2 = \frak Q$ (the ``canonical double cover'') are purely imaginary. This is the so--called {\it Boutroux condition}, \eqref{Boutroux}.
\end{itemize}
See Theorem \ref{quadrdiff}, Corollary \ref{corboutroux}.
This result borrows heavily on ideas from  \cite{BertoGrootKuijlaars}. 

We have also a converse to the above; namely (Theorem \ref{converseBoutroux}) we show that any quadratic differential $\frak Q$ satisfying the same conditions  in the above bullet-points, defines a set of minimal capacity obtained by collecting the graph of some of its  {\it critical vertical trajectories} (for the terminology and background on quadratic differentials we refer to the monographs \cite{StrebelBook, JenkinsBook}).

In Section \ref{secJIP} we prove a  result of a slightly different flavour. Namely we establish a criterion by which the capacity of two sets can be compared to each other. It is immediate to see from the definition and Frostman's theorem (generalized to higher genus) that the capacity is monotonic with respect to inclusion. The comparison criterion we prove is a generalization of a result of Jenkins' \cite{JenkinsArt}, where he proves that sets with the $S$-property\footnote{Jenkins in loc. cit. characterizes the sets in a different way but in the modern terminology his characterization is equivalent to the $S$-property.} are of minimal capacity in an appropriate class.  The generalization consists in extending the result to sets without the $S$-property (we still assume that their boundaries  are sufficiently smooth, however): here we will be using the ideas of \cite{BertoTovbisNLSCritical}.  
This leads to the interesting notion of {\it Jenkins Interception Property} which we outline here. 
Let $\frak F$ be a poly-continuum with piecewise smooth boundary. Its Green's function $G_{\frak F}(p;\infty)$ defines a measure on the foliation of its gradient lines. 
The density of this measure can be expressed as a measure on the outer boundary of $\frak F$ with density given by the normal derivative of the Green's function itself. If, from a smooth point $p$ on the outer boundary, there are two emerging trajectories (from either sides) we call {\it dominant} the trajectory on the side where the above density is larger, the other being the {\it recessive} trajectory. If the density is equal on both sides, as is the case if $\frak F$ has the $S$--property,  we say that both trajectories are dominant. If there is only one trajectory from $p$ we still call it dominant. A second set $\K$ has the Jenkins Interception Property (JIP)  relative to $\frak F$ if it intersects all but finitely many of the dominant trajectories. Then the result is Theorem \ref{thmJIP}, stating that 
\begin{shaded}
if  the poly-continuum $\K$ satisfies JIP relative to $\frak F$ then its capacity is  larger. 
Furthermore, under this assumption, the two capacities are equal iff $\K=\frak F$.
\end{shaded}
This result perfectly illustrates the importance of the $S$-property, when all trajectories are dominant. Indeed, intuitively it clarifies that any deformation of a set with the $S$-property (at fixed anchors) has larger capacity because it intersects all the trajectories either on ``one or the other side''. 

In Corollary \ref{corS} and Proposition \ref{propuniqueness} we use this result to show the uniqueness of the minimizer in a given connectivity class. 

We comment in Section  \ref{CommentJS} on the relationship between the Chebotarev problem and the theory of Jenkins-Strebel differentials, which are used in the combinatorial model of the moduli spaces of Riemann surfaces. 

We conclude in Section \ref{conclusion} outlining in slightly greater detail how the present results fit into a future plan of study of large degree asymptotic behaviour of Pad\'e\ approximation on a Riemann surface.

\paragraph{Acknowledgments.}
The author completed the work during his tenure as Royal Society Wolfson Visiting Fellow (RSWVF\textbackslash R2\textbackslash 242024) at the School of Mathematics in Bristol University. 
The work  was supported in part by the Natural Sciences and Engineering Research Council of Canada (NSERC) grant RGPIN-2023-04747.

\section{Setup: connectivity patterns} 
\label{connpatt}


Let $\CC$ be a compact smooth Riemann surface of genus $g$ with complex structure.
We will also need a complete  metric, $\sigma$, in the same conformal class; any metric will do, e.g. the unique metric of constant negative curvature $-1$. When referring to the distance between points or sets, we will understand it in this metric and we denote the distance function by the symbol $d_{\sigma}:\CC\times \CC \to \R_+$. Corresponding to this metric we have a Hausdorff metric, which we denote by $d_{H}$, on the class of compact subsets of $\CC$. This is also a complete metric. 

Let $E=\{e_1, \dots, e_N\}\subset \CC$ be a finite nonempty subset which we call {\it anchor set} (with the elements being {\it anchors}). 
Let $\pi_1(\CC, E)$ be the relative homotopy group; thus the elements of this group are homotopy classes of paths with endpoints at one or more anchors. 
By {\it poly-continuum} we understand a finite union of {\it continua} (compact connected sets). 
\paragraph{Decorated adjacency matrices.}
In broad terms, we need  to classify the poly-continua containing the anchor set $E$ according to their    ``connectivity pattern", which we now proceed to formalize. 

Let $\gamma\in \pi_1(\CC, E)$;  we say that a poly-continuum $\K$ {\it contains the class $\gamma$} if any $\epsilon$-fattening of $\K$ contains a representative path in the class $\gamma$. It should be clear that the property "continua containing the pre-assigned class $\gamma$" is Hausdorff continuous, namely, if a Hausdorff convergent sequence $(\K_n)_{n\in\N}$ has the property that all $\K_n$'s contain $\gamma$, then so does the limiting continuum $\K_\infty$. 

Given any poly-continuum $\K$ containing the anchor set $E$, we can associate a {\it connectivity pattern}  that records  which points of $E$ belong to the same connected component of $\K$, as was done in \cite{BertoTovbisNLSCritical}. In view of the nontrivial topology of $\CC$, however, we want to retain a more granular information of how two points of $E$ are connected by $\K$. For this reason we introduce the {\it decorated adjacency matrix} below. 

\bd
Given a poly-continuum $\K$ containing the anchor set $E$, the decorated adjacency matrix $M = M(\K)$ is a square matrix of size $|E|\times |E|$, indexed by elements of $E$. In the entry $M_{e,\wt e}$ we place the set of all classes of paths $\gamma\in \pi_1(\CC,E)$  connecting $e$ to $\wt e$ and  $\wt e$ to $e$, and contained in $\K$. 
\ed
We also define a partial order amongst these decorated adjacency matrices as follows; we say that $ M(\wt K)\prec M(\K) $ if the set in each entry of $M(\K) $  contains the set in the corresponding entry of $ M(\wt K)$. In words we will say that {\it the connectivity of $\K$ exceeds that of $\wt K$}. 

Finally we want to define a class of poly-continua that all exceed a given connectivity pattern $\scr P$.

For each pair $(e,\wt e)\in E\times E$ we choose either the empty set or a pair $\{\gamma, \gamma^{-1}\}$ with $\gamma$ a choice of a  unique class $\gamma_{e,\wt e}\in \pi_1(\CC,E)$.
On the diagonal $(e,e)$ there is no constraint and hence the choices are either the empty set, the singleton of the identity or  a loop based at $e$ together with its inverse.  
We construct then a matrix $M^{\scr P}$ containing in each entry the corresponding choice. 

 For our application later on we insist that every anchor $e\in E$ should be the endpoint of a nontrivial class:
\be
\label{cond}
\forall e\in E\ \exists \ \wt e:\  M^{\scr P}_{e,\wt e} \neq \emptyset, {\rm Id}.
\ee

We briefly comment on this condition:  its motivation  is that we are going to consider minimal capacity poly-continua $\K$ such that $M^{\scr P}\!\!  \prec M(\K)$, and if an element $e$ does not belong to a nontrivial class in the pattern $\scr P$, then  this point would be ``ignored'' by the minimal capacity set (isolated points have zero capacity), and we may as well have removed $e$ altogether from consideration. 

We thus come to the formal definition of our class of poly-continua:
\begin{definition}
\label{defclass}
Let $\scr P$ be an admissible pattern as per condition \eqref{cond}. 
Let $\mathbb K_{\scr P}$  be the class  of all poly-continua $\K$ such that: 
\begin{enumerate}
\item the connectivity pattern of $\K$ exceeds $\scr P$, namely, $M^{\scr P} \prec M(\K)$;
\item the number of connected components of $\K$ does not exceed the number of homotopy classes in $\scr P$ modulo inversion (hence at most $|E|(|E|+1)/2$).\hfill $\blacktriangle$
\end{enumerate}
\end{definition}
Again, a comment is in order here. While the first condition in Definition \ref{defclass} should be self-explanatory, the second condition is a minimality of sorts. Namely, we insist that all components of $\K$ are ``actively participating" in connecting anchors. This is convenient because in the minimization process of the capacity, we want to avoid the occurrence of components of small capacity.

The class $\mathbb K_{\scr P}$ is closed in Hausdorff topology since we have already commented that the condition of containing a class is a Hausdorff continuous property and so is the uniform bound on the number of connected components.

\begin{problem}
\label{problem!}
For a given anchor set $E$ and admissible connectivity pattern $\scr P$ we are seeking the set of minimal capacity in $\mathbb K_{\scr P}$.
\end{problem}

For Problem \ref{problem!} to make sense we need do carefully define the notion of capacity on a Riemann surface, and we take most of what follows from \cite{Chirka}.

\section{Capacity on a Riemann surface}
\label{secCapacity} 
As in the standard case where $\CC$ is the Riemann sphere, the capacity of a set is defined relative to a chosen point in its complement. In the case of the Riemann sphere this is usually taken as the point at infinity, but for a Riemann surface $\CC$ there are no distinguished points. We thus choose a point  on $\CC$ which we denote $\infty$, where the notation is simply a nod to the standard case.

In all considerations in this paper we will assume $\K$ to be a poly-continuum with finitely many connected components. It is easy to see that such a set satisfies the Wiener criterion for the well--posedness of the Dirichled problem (see \cite{SaffTotik}, App A.2, Theorem 2.1).
\begin{definition}
\label{defGK}
Given a proper poly-continuum $\K\subset \CC$ with $\infty\not\in\K$,  the polar Green's function of $\K$  (with pole at $\infty$) is the unique function $G_\K(p;\infty)$ such that 
\begin{enumerate}
\item $G_\K(p;\infty)$ is harmonic in the variable $p\in \K\setminus \{\infty\}$;
\item $ G_\K(p;\infty) = 0 $ for all $p\in \K$;
\item $G_\K(p;\infty) - \ln \frac 1 {|z_\infty(p)|}$ is harmonic and bounded in a full neighbourhood of $\infty$, where $z_\infty$ is any local coordinate near $\infty$. \hfill $\blacktriangle$
\end{enumerate}
\end{definition}
By the minimum principle $G_\K(p;\infty)$ is a strictly positive function on $\CC \setminus \K\cup\{\infty\}$.

Observe that $G_\K$ is non-trivial only in the connected component,  $\Omega$,  of $\CC \setminus \K$ that contains $p=\infty$. We call this the {\it exterior region} and we call $\pa \Omega$ the {\it outer boundary} of $\K$.  
Different compact sets  with the same $\Omega$ will have the same Green's function, namely, the Dirichlet problem is ``blind'' towards holes of $\K$. We will say that $\K$ is {\it polynomially convex} if $\K= \CC\setminus \Omega$, although the  adverb ``polynomially'' is here  a bit out of place, since there are no polynomials on $\CC$. To simplify the discussion we make the running assumption that 
\be
\label{ass}
\text{all poly-continua we  consider will be  polynomially convex.}
\ee

\subsection{Bipolar Green's function}
We need a substitute for the Green's kernel $\ln | z-w|$ on the Riemann sphere. Following \cite{BertoGrootKuijlaars} we introduce  the notion  below.
Let $G_\infty(p,q)$ be the unique (up to additive constant) {\it bipolar Green function} with logarithmic pole at $\infty$. Namely 
\begin{enumerate}
\item $G_\infty(p,q)$ is harmonic w.r.t. $p$ in $\CC\setminus \{\infty, q\}$;
\item in any  local coordinate, for $p$ near $q$
\be
G_\infty(p,q)= \ln {|z(p)-z(q)|} + \mathcal O(1).
\ee
Equivalently, $G_\infty(p,q)-\ln d_\sigma(p,q)$ is bounded along the diagonal $p=q\neq \infty$.
\item in a local coordinate near $\infty$,  $z_\infty(\infty)=0$, we have 
\be
G_\infty(p,q) =\ln \frac 1  {|z_\infty(p)|} + \mathcal O(1).
\ee
Equivalently, $G_\infty(p,q)+\ln d_\sigma(p,\infty)$ is bounded for $p$ in a neighbourhood of $\infty$. 
\end{enumerate}
The bipolar Green's function can be constructed as follows; let $\Omega^{i\R}_{\infty, q}$ be the unique differential of the third kind  with purely imaginary periods and with simple poles at $p= q, \infty$ and residues $1, -1$,  respectively \cite{Fay}.   Then 
\be
\label{GOmega}
G_{\infty}(p,q) =\Re \le( \int^p \Omega^{i\R}_{\infty,q}(\zeta)\ri). 
\ee
up to additive constant.
Denote with  $\pa_p$  the Wirtinger operator that takes any local harmonic function $F(p)$ and returns a local holomorphic differential given by 
\be
\label{Wirtinger}
\pa_p F(p) := \frac 1 2 \le(\frac{\pa F}{\pa x} - i \frac {\pa F}{\pa y}\ri) \d z, \ \ z = x+iy,
\ee
in a local coordinate $z=z(p)$. Then we have the converse of \eqref{GOmega} 
\be
\label{OmegaG}
\Omega^{i\R}_{\infty,q} (p) = 2\pa_p G_\infty(p,q).
\ee
It  can also be shown either directly or by a consequence of the Riemann bilinear identities that $G_{\infty}(p,q) = G_\infty(q,p)$. 
The arbitrariness of additive constant can be disposed of in several ways. Possibly the most expedient is to simply give an explicit formula for $G_\infty(p,q)$ in terms of Riemann Theta functions:
{ \be
\label{expG}
G_\infty(p,q) = \ln \le|\frac 
{\Theta_\Delta\le(\int_q^p\bs\omega \ri)}{ \Theta_\Delta\le(\int_\infty^p\bs \omega \ri) \Theta_\Delta\le(\int^\infty_q\bs \omega \ri)}\ri| -2 \pi
\Im \le(\int_\infty^p\bs \omega^t\ri) \cdot (\Im{\bs\tau})^{-1} \cdot \Im \le(\int_\infty^q\bs \omega\ri)
\ee}
where $\bs \omega = [\omega_1,\dots, \omega_g]$ is the vector of normalized holomorphic  differentials and  $\Theta_\Delta$ denotes the Riemann theta function with an odd-nonsingular characteristic, see \cite{Fay}, Ch. I and also Appendix \ref{expCauchy} here. 
\begin{example}
If $\CC$ is an elliptic curve (genus one) represented by $\C /\Z+\tau \Z$ with $\Im \tau>0$ and with the choice of $\infty$ as $p=0$, then (see also \cite{BertoGrootKuijlaars} Example 2.3)
\be
G_\infty(p,q) =  \ln \le|\frac 
{\theta_1\le(p-q\ri)}{ \theta_1\le(q \ri) \theta_1\le(p\ri)}\ri| -2 \pi
\frac {\Im \le(p\ri) \Im \le(q\ri)}{\Im{\tau} },
\ee
where $\theta_1$ is the standard Jacobi elliptic theta function. 
\end{example}
 Alternatively one can choose the additive constant subordinate to the choice of a local coordinate $z_\infty$ and so that $G_\infty(p,q) = \ln  \frac 1{|z_\infty(p)|} + \mathcal O(|z_\infty(p)|)$ (namely, so that the constant term vanishes in the coordinate $z_\infty$).   
\subsection{Capacity of a  set $\K$}
Recall the running assumption that $\K$ is a poly-continuum with finitely many components. 
\bd
\label{defGreenK}
Given the Green's function $G_\K(p;\infty)$ and a fixed local coordinate $z_\infty$ we define  the capacity as 
\be
\ln {\rm Cap}_{z_\infty}(\K):= \lim_{p\to \infty} \le(\ln \frac 1{|z_\infty(p)|} -G_\K(p;\infty) \ri). 
\ee
Alternatively and equivalently, in the local coordinate $z_\infty$ the Green's  function has the expansion
\be
\label{GCap}
G_\K(p;\infty) = \ln \frac 1{|z_\infty(p)|} - \ln {\rm Cap}_{z_\infty}(\K) + \mathcal O(|z_\infty(p)|). 
\ee
\ed
One could tie the choice of coordinate to the metric by choosing the local coordinate adapted to the metric in the sense that $\lim_{p\to\infty}\frac {|z_\infty(p)|}{d_\sigma(p,\infty)} = 1$, but there is no discernible convenience in doing so.

It is useful to give several equivalent representations of this capacity and relate it to a bipolar Green's function. 
%
%


For the case $\CC = \mathbb P^1$ we have the classical Frostman's theorem \cite{Frostman} relating the Green's function of  the complement of a compact set   to a probability measure on its boundary. 
In particular, such a measure is also the unique minimizer of the logarithmic energy of probability measures supported on the compact set. 

This has been generalized in \cite{Chirka} to an arbitrary Riemann surface with even greater generality than needed in our paper now. 
We briefly review the key points, with the intention of simplifying the reader's navigation.

In our setting, it follows  that there is a probability measure $\d\mu_\K $ supported on the  outer boundary of $\K$ 
\be
\label{GGG}
G_\K(p;\infty) =\mathcal U_{\d\mu_{_\K}}(p) - \ln {\rm Cap}_{z_\infty}(\K), \qquad
\mathcal U_{\d\mu}(p):=  \int  G_\infty(p,q)\d \mu (q) 
\ee
where the bipolar Green's function is normalized so that $G_\infty(p,q) = \ln \frac 1 {|z_\infty(p)|}+\mathcal O(|z_\infty(p)|)$ and  the constant ${\rm Cap}_{z_\infty}(\K)$ is precisely the previously defined capacity subordinated to the same choice of local coordinate. 
Integrating \eqref{GGG} against $ \d \mu_\K $ and recalling that by definition $G_\K$ vanishes on $\K$, we obtain the expression of ${\rm Cap}_{z_\infty}(\K)$ as follows
\be
\ln {\rm Cap}_{z_\infty}(\K )=  \iint G_\infty(p,q) \d \mu_\K (p) \d \mu_\K (q) .\label{GreenI}
\ee
The analog of Frostman's theorem, following \cite{Chirka} and adapting to our situation here, is summarized in the following proposition.

\bp
\label{propFrostman}
Let $\M_\K^1$ denote the set of probability measures supported on $\K$. Then 
\bea
 {\rm Cap}_{z_\infty}(\K) &= {\rm exp} \le(-\ds \mathcal I(\K)\ri),\ \ \ 
 \mathcal I(\K):= \min_{\d \mu\in \M_\K^1}\mathcal E[\d \mu]\nn\\   
&\mathcal E[\d \mu]:= -\int_{\K\times \K} G_\infty(p,q) \d \mu(p) \d \mu(q)
\eea
The minimum of $\mathcal E$ is achieved at a uniquely defined  measure $\d \mu_\K $, supported on its outer boundary. The Green's function is related to the  potential of the measure by \eqref{GGG}, i.e., they differ by a constant equal to  the logarithm of the capacity.
\ep
In the classical case, as mentioned, this is essentially Frostman's theorem; in the case of Riemann surfaces instead, with reference to \cite{Chirka} and notation therein:
\begin{enumerate}
\item  The existence and uniqueness of $\d\mu_{\K}$ can be found in  Section 2. 
In order to adapt loc.cit. to our case, we should take $\phi = \delta_\infty$. See Lemmas 1,2 in section 2.3 in particular.
\item The fact that the measure is supported on the outer boundary of $\K$ (in our terminology) and that the potential  is equal to the Green's function (up to constant) follows then from  Corollary 3, pag 319. 
\end{enumerate}
The quantity $\mathcal E[\d \mu]$ can be termed the (electrostatic) {\it energy} of a probability measure, in analogy with the case of the Riemann sphere where the bipolar Green's function would be $\ln  {|z-w|}$. %
\subsection {The conformal map $W_\K$.}
\label{secunif}
Let $g_\K(p;\infty):= G_\K(p;\infty) + i H_\K(p;\infty)$, where $H_\K$ is the harmonic conjugate function so that $g_\K$ is analytic and defined on the universal cover of $\CC \setminus \K\cup \{\infty\}$. It has a logarithmic branch point at $\infty$ and with additive monodromy $2i\pi$, behaving like $\ln  {\z_\infty(p)}$, where $\z_\infty(\infty)=0$ is a local parameter. 
Since $\Re g_\K = G_\K$ is zero on $\K$ the function 
\be
\label{defW}
W(p) = W_\K(p) := {\rm e}^{-g_{_\K}(p)}
\ee
maps the universal cover of $\CC \setminus \K$ onto the unit disk, with $\infty$ being mapped to the origin in the $W$--plane and the boundary of $\K$ onto the unit circle. In terms of a local coordinate $z_\infty$ near $\infty$ we have 
\be
W(p) = {z_\infty(p)}{{\rm Cap}_{z_\infty}(\K)}  \Big(1+ \mathcal O(z_\infty)\Big) =  {z_\infty(p)}{\rm e}^{-\I(\K)}  \Big(1+ \mathcal O(z_\infty)\Big).
\ee
\br
When $\CC=\mathbb P^1$ (the Riemann sphere), the bipolar Green's function is simply $\ln { |z-w|}$ and the subtleties related to the choices of coordinate near $\infty$ are lost. When referring to the logarithmic energy, it is universally implied that the natural choice of coordinate $z$ has been made near $\infty$.
\er

Consider the map $W_\K$ defined in \eqref{defW}. In general it is not one-to-one unless we perform some additional dissection that we now describe.
 The Green's function $G_\K(p;\infty)$ must have a certain number $k\geq 0$ of  critical points (counted with multiplicity); these are the zeros of the meromorphic differential $\pa_p G_\K(p;\infty)$, where $\pa_p$ is the Wirtinger operator \eqref{Wirtinger}. This differential is meromorphic on $\CC\setminus \K$ with a unique simple pole at $\infty$ and residue $-\frac 1 2$.  The number of these zeros (counted with multiplicity) depends on the Euler characteristic of $\CC \setminus\K \cup\{\infty\}$ thanks to elementary Morse theory. 
Namely, the degree of the zero divisor of $\pa_p G_\K $ on $\CC\setminus \K$ is  
\be
 k = - \mathcal X(\CC \setminus \K\cup\{\infty\}). 
\ee
\begin{example}
If $\K$ is a single contractible component  then $k = 2g$. 
\end{example}
The map $W = W_\K$ determines the flat metric by the formula  $ |2\pa_p G_\K| = \le|\frac {\d W}{W}\ri|$ on $\CC \setminus \K\cup\{\infty\}$, with conical singularities at the $k$ critical points $\scr B = \{b_1,\dots, b_k\}$. 

In this metric both the level-sets  and the gradient lines of $G_\K$ provide  two  local foliations in mutually orthogonal geodesics. The dissection, $\Sigma$, is then constructed as follows;
from each critical point $b\in \scr B$ there are $2\mu_b$ steepest descent gradient lines of $G_\K$, where $\mu_b$ is the multiplicity of the zero of $\pa_pG_\K  $ at $b$. Each of these gradient lines ends either on another critical points, or on the boundary of $\K$.

We define $\Sigma$ as the union of all these steepest descent trajectories from all critical points in $\scr B$. 

\bl
\label{lemmastar}
The set $\CC \setminus( \Sigma \cup \K)$ is simply connected. 
The map $W_\K$ maps $ \CC \setminus \Sigma \cup \K$ conformally and bijectively to the unit disk with $2k$  radial slits. 
\el
\noindent{\bf Proof.}
The  proposed set is geodesically  star-like with respect to the flat metric $|\d W|/|W|$, with the point $\infty$ at the center of the star; for each $p_0\not\in \Sigma\cup \K$ we take the maximally extended steepest ascent gradient line, $\gamma$,  of $G_\K$. Any steepest ascent trajectory of $G_\K$ can end only at $\infty$ or at a critical point, but since we have precisely removed all the steepest ascent geodesics ending at critical points, the geodesic $\gamma$ above  can only end at $\infty$, and hence  it is mapped to a radial segment in the $W$--plane, ending at the origin. This proves the simple connectedness. 

The second assertion now follows from the fact that $ W_K$ restricted to $\CC\setminus \Sigma \cup \K$ is an analytic map on a simply connected domain without any critical points, hence bijective. 
\QED
We point out  that if $\CC \setminus \K$ is already simply connected, then $W_\K$ will be the corresponding  uniformizing map mapping $\infty$ to $W_\K=0$. 
\section{Generalized Chebotarev problem in higher genus}
\label{secChebo}
In this section we prove the first main result of the paper, Theorem \ref{thmP}.
We consider an admissible  connectivity pattern $\scr P$ and the corresponding class $\mathbb K_{\scr P}$ as in Definition \ref{defclass}. 
We fix a local coordinate $z_\infty$ and consider then the corresponding capacity as a functional 
\be
{\rm Cap}_{z_\infty} : \mathbb K_{\scr P} \to \R_+.
\ee
We will  show that there exists a minimizer of the capacity within this class. We then investigate the properties of such minimal sets in the subsequent sections, including the uniqueness of the minimizer. 

This problem can be considered as a generalized Chebotarev problem in higher genus. Let us recall that the classical Chebotarev problem consists in  finding a continuum $\K\subset \C$ containing a finite set $E$ (our ``anchors'') and minimizing the usual logarithmic capacity \cite{Polya}.  By  {\it generalized} Chebotarev problem (in genus zero) we mean that we allow poly--continua exceeding a given connectivity pattern, everything else being equal, see also \cite{MFRakh}.   

When formulating the problem in higher genus we must first define the relevant notion of capacity, which we have done in the previous section.  In addition we need to contend with the nontrivial topology of $\CC$, which is the reason behind  the rather convoluted description of connectivity patterns.

We need two preliminary lemmas. 
For any $\delta>0$ (sufficiently small)  let 
\be
\mathbb K_{\scr P} ^{(\delta)} := \{ \K\in \mathbb K_{\scr P} :\ \ \ d_\sigma(\K,\infty)\geq \delta\},
\ee
namely, those poly-continua that exceed the connectivity pattern $\scr P$ and stay away from $\infty$.  
The following Lemma \ref{lemmacont} is  a straightforward consequence of the results in \cite{Chirka} (in slightly different terms). See the (unnumbered)  Lemma on page 332 therein. To help the reader navigate the reference, the surface is denoted $S$,  the  probability measure $\phi$ in the Lemma should be taken as the delta function supported at $\infty$. See also Corollary 1 on page 333, where the family denoted $\mathcal X$ there should be taken as our $\mathbb K_{\scr P}^{(\delta)}$. 
\bl
\label{lemmacont}
Let $\delta>0$ be sufficiently small so that
 $\mathbb K_{\scr P}^{(\delta)} $ is not empty. Then the capacity 
 $
 {\rm Cap}_{z_\infty} : \mathbb K_{\scr P} ^{(\delta)}\to \R_+
 $ is Hausdorff continuous. 
\el
We sketch an elementary proof in Appendix \ref{AppProofcont}, for the interested reader. 
The next lemma is also essentially in \cite{Chirka} (last bullet point on page 307) but we provide a slightly different proof nevertheless.
\bl 
\label{lemmainfty}
Suppose that $\K_n$ is a  sequence of poly-continua in $\mathbb K_{\scr P}$ such that $d_\sigma(\K_n,\infty)\to 0$. Then 
$$
\lim _{n\to\infty} {\rm Cap}_{z_\infty}(\K_n)=+\infty.
$$  
\el
\noindent {\bf Proof.}
 We transfer the problem to one of the standard complex plane. Fix $r>0$ sufficiently small so that   the set $|z_\infty|\leq r$ in $\CC$  is mapped conformally to the disk $\DD_r$ and so that all anchors are outside.  The map $\zeta = \frac 1{z_\infty}$ then maps  to  $|\zeta|\geq \frac 1 r$. 
In the coordinate $\zeta$ the bipolar Green function   can be written 
\be
-G_\infty(p,q) = \ln \frac 1{|\zeta - \xi|} + H(\zeta,\xi),
\ee
where $\zeta = \frac 1{z_\infty(p)}, \xi = \frac 1 {z_\infty(q)}$ and $H(\zeta,\xi)$ is a function which is harmonic and bounded in each variable  in $|\zeta|\geq \frac 1 r \leq |\xi|$. Let us set 
$$  M:= \sup _{|\zeta|, |\xi|>R} \!\!\! H(\zeta,\xi)<\infty,\ \ \ R:= \frac 1 r.$$
We are going to construct a sequence of probability measures supported within 
$$\K_n^r:= \K_n\cap \{|z_\infty|\leq r\}$$  whose energy tends to $-\infty$.   The map $|\zeta|$ is continuous and hence maps each connected component of  $\K^r_n$ to an interval. Furthermore, since each such connected component must contain also an anchor, it must ``exit'' the disk $|z_\infty|<r$. Thus, each component of $\K_n^r$ is mapped to an interval $[R,c_n]$. Then $[R,L_n]$ is simply the largest of the (finitely many!) such intervals, with $L_n = \frac 1 {\inf _{\K^r_n}|z_\infty(p)|}$. Note that $L_n$ tends to infinity since the distance of $\K_n$ from $\infty$ tends to zero. 

Let $\d \rho_n$ be the measures obtained  by (any) lifting  to $\K_n^r$ by the map $x = |\zeta|$ of the arcsin measure $\frac { \d x} {\pi\sqrt {|(x-R)(x-L_n)|}}$ on $\R_+$ restricted to the intervals $[R,L_n]$. 
Existence of a lifting is a standard fact,  see for example the lecture notes  \cite{Notes}, sec. 14.3.   By the standard inequality $
|\zeta - \xi| \geq\Big||\zeta | - |\xi| \Big|$ we have  that  $\ln \frac 1{|\zeta - \xi|} \leq \ln \frac 1{\Big||\zeta | - |\xi| \Big|}$. Now, 
\bea
\nn
- \ln {\rm Cap}(\K_n) &= \inf _{\d \rho\in \M^1_{\K_n}} -\iint \d\rho(p)\d\rho(q) G_\infty(p,q)\leq 
\\ \nn
&\leq \iint \d\rho_n(\zeta) \d \rho_n(\xi) \ln \frac 1{|\zeta - \xi|} + \iint \d\rho_n(\zeta) \d \rho_n(\xi) H(\zeta, \xi) \leq 
\\
&\leq M+ \int_{R}^{L_n} \ln \frac 1{|s - t |} \frac {\d s}{\pi\sqrt{|s-R||s-L_n|}}\frac {\d t}{\pi\sqrt{|t-R||t-L_n|}} = M -\ln \frac {L_n-R} 4.
\eea
The last equality is due to the well known fact that the arcsin measure realizes the minimal energy of the segment, together with the fact that the capacity of a segment is its length divided by four. Since $L_n\to \infty$ the capacity must diverge. 
\QED

With the two lemmas established, we can proceed to formulate and prove the main theorem of this section.

\bt
\label{thmP}
For any admissible connectivity pattern $\scr P$ there exists a (unique) continuum $\K_{\scr P}$ minimizing the  capacity within $\mathbb K_{\scr P}$. 
\et
\noindent {\bf Proof.}
Let  
$
C_{\scr P}= \inf {\rm Cap}_{z_\infty}(\mathbb K_{\scr P}) 
$ and 
let $(\K_n)_{n\in\N} \subset \mathbb K_{\scr P}$ be any  sequence of poly-continua such that ${\rm Cap}_{z_\infty}(\K_n) \to C_{\scr P}$. 
By Lemma \ref{lemmainfty} any such a sequence must eventually maintain a positive distance  $\delta$ from $\infty$ and  hence  belongs to $\mathbb K_{\scr P}^{(\delta)}$. Since such a class is Hausdorff compact and ${\rm Cap}_{z_\infty}$ is continuous  on it, it must attain a minimum. 
We are not proving here the uniqueness. This will be established in Section \ref{secJIP}, Proposition \ref{propuniqueness}. 
\QED

\section{Properties of sets of minimal capacity}
\label{secminimal}
The next question is to characterize the geometry of these minimizers; for the case $\CC=\mathbb P^1$ it is well known that these consist of finite union of analytic arcs and that they are related to the so--called horizontal trajectories of quadratic differentials satisfying the so-called $S$-property,  see Stahl \cite{Stahl, Stahl2, Stahl3} or \cite{MFRakh}. 

A completely similar situation occurs here as we are discussing in this section

The key tool is, as typical, the notion of {\it Schiffer variation} \cite {SS54}, which was used in higher genus already in \cite{BertoGrootKuijlaars}. 
We recall this notion here.

Let $h$ denote a smooth vector field on $\CC$; it determines a group of diffeomorphisms 
$$
\Phi: \R\times \CC \to \CC
$$
 with the property that $\frac {\d \Phi}{\d t}(t,p)\bigg|_{t=0}  =h(p) $ and $\Phi(0,p) = p$.  

For any measure $\d\mu$ one then pulls back the action of this diffeomorphism to define a transported measure $\d\mu_{\epsilon, h}$ according to the formula
\be
\int f(p) \d\mu_{\epsilon,h}(p) := \int f(\Phi(\epsilon,p)) \d\mu(p), \ \ \epsilon \in \R. 
\ee 
The following is the adaptation of loc. cit.
\bd[Definition 3.3 of \cite{BertoGrootKuijlaars}]
A mesure $\d\mu$ on $\CC$ with support not containing $\infty$ is a {\bf critical measure} if 
\be
\lim_{\epsilon \to 0}\frac { \mathcal E[\d\mu_{\epsilon, h}] - \mathcal E[\d\mu]}\epsilon =0
\ee
for all vector fields $h$ at least $C^1$ in a neighbourhood of the support of $\mu$.  Here $\mathcal E[\d\mu]$ is the energy functional (see also Prop. \ref{propFrostman})
\be
 \mathcal E[\d \mu]:= -\iint  G_\infty(p,q) \d \mu(p) \d \mu(q).
 \label{energy}
\ee 
\ed
The following proposition is given without proof because the proof is verbatim the one in loc. cit.
\bp [Prop. 3.4 in \cite{BertoGrootKuijlaars}]
\label{propvariation}
Let $\d\mu$ be a measure of finite energy $\mathcal E[\d\mu]<+\infty$ \eqref{energy}. Then 
\begin{enumerate}
\item For any vector field $h$ which is $C^1$ in a neighbourhood of the support of $\d\mu$ the limit below exists
\be
\lim_{\epsilon \to 0}\frac { \mathcal E[\d\mu_{\epsilon, h}] - \mathcal E[\d\mu]}\epsilon  = \Re D_{h}(\d\mu),
\ee
where $D_h(\d\mu)$ is the Schiffer variational derivative along the flow generated by the vector field $h$ and it is given by 
\be
\label{critcond}
D_{h}(\d\mu) = \iint h(p)\Big( \pa_p G(p,q) + h(q) \pa_q G(q,p)\Big)\d\mu(p)\d\mu(q).
\ee
\item The measure is critical if $D_h(\d\mu)=0$ for any vector field $h$ which is $C^1$ in a neighbourhood of the support. 
\end{enumerate}
\ep
With reference to Remark 3.6 in loc. cit. we recall that it suffices to verify the criticality under variations induced by vector fields that are in the holomorphic tangent bundle to the surface, namely, for vector fields that in a local coordinate can be written as 
\be
h = g(z, \ov z) \frac {\pa}{\pa z}.\nn
\ee
We also recall that given any locally holomorphic differential $\omega = f(z) \d z$, the expression $\frac1  \omega$ can be interpreted as a possibly meromorphic vector field in the holomorphic tangent bundle.  In fact we will use only {\it meromorphic} vector fields, which are thus simply the reciprocal of a (meromorphic) differential.

The idea to be implemented is that by choosing the vector field appropriately, we can derive properties of the differential of the potential of a critical measure.

\paragraph{Preparation.}
In order to provide an appropriate supply of meromorphic vector fields over which to test the criticality condition \eqref{critcond}, we need a special kernel which is denoted $C^{(2,-1)}(p,q)$ also in loc. cit. but was only defined for the case of an elliptic Riemann surface,  and now has the properties described below.
The main difference and novelty relative to \cite{BertoGrootKuijlaars} is that the transporting diffeomorphism must preserve the anchor set $E$, namely
\be
\text{the criticality condition must hold for all $h$ that vanish at the anchor set $E$.}
\ee
This is the same constraint as in \cite{MFRakh}. 

Let $\scr B_{3g}$ denote a generic divisor of degree $3g$. Then 
\begin{enumerate}
\item $C^{(2,-1)}(p,q)$ with respect to $p$ is a quadratic differential such that 
\begin{enumerate}
\item has simple pole for $p=q$ and at $E$;
\item has a zero of order $N-3$ at $\infty$ and $3g$ other zeros at $\scr B_{3g}$.
\end{enumerate}
\item With respect to $q$ it is a $-1$ differential (a meromorphic vector field) such that 
\begin{enumerate}
\item has a pole a $q=p$, at $q\in \scr B_{3g}$ and a pole of order $N-3$ at $q=\infty$;
\item has zeros at $E$.
\end{enumerate}
\item it is normalized in such a way that, in any local coordinate $\zeta$
\be
C^{(2,-1)}(p,q) = \frac {\d z^2}{\d w} \le (\frac 1{z-w} + \mathcal O(1)\ri), \ \ z\to w,
\ee
where $z=\zeta(p),\ w=\zeta(q)$.  This normalization is coordinate independent.
\end{enumerate}
Using the common notation for divisors we can summarize the above in the formula:
\bea
\div_ p \le(C^{(2,-1)}(p,q)\ri)\geq -(q) - E + (N-3)(\infty) + \scr B_{3g}, \nn\\
\div_q \le(C^{(2,-1)}(p,q)\ri)\geq -(p) + E - (N-3)(\infty) - \scr B_{3g}.
\eea
The existence and uniqueness, under genericity assumption,  is guaranteed by the Riemann--Roch theorem. An explicit formula for $C^{(2,-1)}$ can be given but  is not essential. See Appendix \ref{expCauchy}.

Consider the expression:
\be
\label{defF}
F(x;p,q):= -C^{(2,-1)}(x,p) \Omega^{i\R}_\infty(p,q)- C^{(2,-1)}(x,q)\Omega^{i\R}_\infty(q,p) + \Omega^{i\R}_\infty(x,p) \Omega^{i\R}_\infty(x,q),
\ee
where $\Omega^{i\R}$ is the third kind normalized differential with purely imaginary periods, related to $G_\infty(p,q)$ by \eqref{OmegaG}.
It is manifestly a quadratic differential with respect to $x$ and has the following analytic structure. 
\begin{enumerate}
\item It is regular at $x=p$  and at $x=q$. This is verified easily in a local coordinate.
\item It has simple poles at $E$ and a double pole at $x=\infty$ with unit bi-residue and  coming only from the last term.
\item It is a harmonic function with respect to $p$ and $q$ on $\CC\setminus \{\infty\}$
\end{enumerate}

We then formulate the next main result
\bt
\label{quadrdiff}
Let $\scr P$ be any admissible connectivity pattern for the set of anchors $E$; let $\K = \K_{\scr P}$ be the set minimizing the capacity within the class $\mathbb K_{\scr P}$ as per Theorem \ref{thmP}. Let $\d\mu$ denote the equilibrium measure on $\K_{\scr P}$. 
Let $G_\K$ denote the corresponding Green's function as per Definition \ref{defGreenK}.
Then 
\be
\le(2\pa_x G_\K(x;\infty)\ri)^2 = \mathfrak Q(x)
\ee
is a meromorphic quadratic differential on the whole $\CC$ with at most simple poles at $E$, a double pole at $\infty$ with bi-residue $1$. The quadratic differential $\frak Q$ is given by 
\be
 \frak Q (x)= \iint F(x;p,q) \d\mu(p)\d\mu(q), 
\ee
with $F$ as in \eqref{defF}.
\et
\noindent
{\bf Proof.}
The measure $\d\mu$ is critical and hence we apply the criterion \eqref{critcond} using 
$$
h(p)=  C^{(2,-1)}(x,p).
$$
Here the divisor $\scr B = \scr B_{3g}$ of degree $3g$  used to define the kernel $C^{(2,-1)}$ can be chosen arbitrarily provided that it does not intersect the set $\K$. 
This yields the identity
\be
\label{main1}
0 = \iint \le[C^{(2,-1)}(x,p)\Omega^{i\R}_\infty(p,q) + C^{(2,-1)}(x,q)\Omega^{i\R}_\infty(q,p)  \ri] \d\mu(p) \d\mu(q).
\ee
In the integrand we now use \eqref{defF} and obtain 
\be
\Big(2\pa_x G_\K(x;\infty) \Big)^2=\le(\int\Omega^{i\R}_\infty(x,q) \d\mu(q)\ri)^2 = \iint F(x;p,q) \d\mu(p)\d\mu(q).
\ee
We have established that $F(x;p,q)$ is a meromorphic quadratic differential in $x$ with poles at $E$ and double pole at $\infty$ with unit bi-residue.  Thus the theorem is established with  $\mathfrak Q(x)= \iint F(x;p,q) \d\mu(p)\d\mu(q)$.
In principle the formula holds for $x$ outside of the support of $\d\mu$. However, using the same arguments as in \cite{MFRakh}, Lemma 5.2, we can conclude that it holds everywhere.  
\QED

The theorem, exactly as in the case of minimal capacity  sets on the complex plane, implies that $\K= \K_{\scr P}$ consists of union of vertical critical trajectories of $\mathfrak Q$ see Corollary \ref{corboutroux}; for background on quadratic differentials and their trajectories see \cite{StrebelBook, JenkinsBook}. 
As is common in Teichm\"uller theory, to any  meromorphic quadratic differential $\mathfrak Q$ we can associate a Riemann surface $\wh \CC$ obtained as the locus within the holormophic co-tangent bundle $T^*\CC$ described by the equation 
\be
\wh \CC:= \{ v^2 = \mathfrak Q, \ \ v\in  T^*\CC\} \subset T^*\CC.
\label{doublecover}
\ee
Let us denote by 
\be
\pi:\wh \CC \to \CC
\ee
the canonical projection from the double cover to the base curve. We adopt the convention that $p,q,$ etc.,  denote points of $\CC$ and $\wh p,\wh q, $ etc.  denote points of $\wh \CC$ projecting onto those denoted by the corresponding letter: $\pi(\wh p) = p$. 
On the double cover $\wh \CC$ there are two points $\wh \infty_\pm$ corresponding to $\infty\in \CC$; they are distinguished by the properties
\be
\res_{\wh \infty_\pm} v= \pm 1, \ \ \ \ v:= \sqrt{\mathfrak Q}.
\ee 
For completeness of information we recall the standard construction of $\wh \CC$ by explicitly defining the local coordinates near the zeros and poles of $\frak Q$. Namely, if $p$ is a zero  or pole  of $\mathfrak Q$ of order  $\mu\in \Z $ then, in a local coordinate $z$,  
$$
{\mathfrak Q} = z^{\mu}\le(C^2 + \mathcal O(z)\ri) \d z^2, \ \ \  C \neq 0.
$$
Correspondingly we define a local coordinate $\zeta$ on $\wh \CC$ by 
\be
\label{zetaz}
\zeta^{2\mu+2 } = \le(\frac {\mu/2+1}{C}\int_0^z \sqrt{\mathfrak Q}\ri)^2 = {z^{ {\mu+2}}}\le(1 + \mathcal O(z)\ri),\ \ \  \mu\geq -1.
\ee
If $\mu=-2$ and $\frak Q = \frac {c^2}{z^2} (1+ \mathcal O(z)) \d z^2$ and we choose the determination so that $\sqrt{Q} = \pm \frac c z(1 + \mathcal (z))\d z$ then a coordinate is defined by 
\be
\zeta = {\rm e}^{\pm \frac 1 c \int^z \sqrt{\frak Q}}, 
\ee
where the constant of integration is not important. 
If $\mu<-2$ (not a situation we are encountering here, but we add it for completeness), then we use the same \eqref{zetaz} but integrating from a point in the neighbourhood of zero. 

\br
Although it is a  subtlety that may not be very apparent, in this definition  only the points with $\mu\in 2\Z+1$ are branch-points and lift to a single point on $\wh \CC$, while  all points with  $\mu\in 2\Z$  correspond to two points on $\wh \CC$. This means that for an even-multiplicity zero we are actually {\it resolving} the singularity of the equation $v^2= \frak Q$. So the definition above is actually the curve obtained by resolving the equation $v^2=\frak Q$ in the cotangent bundle of $\CC$, namely, the {\it normalization} of the algebraic set. 
\er
We then have the following corollary.
\bc
\label{corboutroux}
Let $\K = \K_{\scr P}\in \mathbb K_{\scr P}$ be a minimizing set for the connectivity pattern $\scr P$. 
Let $\mathfrak Q$ be the quadratic differential of Theorem \ref{quadrdiff}.
\begin{enumerate}
\item 
The set $\K$ consists of a union of finitely many analytic arcs that coincide with some of the critical vertical trajectories of $\mathfrak Q$ and containing all anchors, where $\mathfrak Q$ has at most a simple pole. 
\item 
The Green's function $G_\K(p;\infty)$ extends harmonically to a function $\Phi(\wh p)$ to the double cover $\wh \CC \setminus \{\infty_\pm\}$, with the property that 
\be
\lim_{\wh p\to\wh \infty_\pm} \Phi(\wh p) = \pm \infty.
\ee
\item 
If $\wh p \mapsto \wh p^\star$ denotes the holomorphic involution of $\wh \CC$ that swaps the sheets, then 
\be
\Phi(\wh p^\star) = -\Phi(\wh p ).
\ee
\item 
The Green function is $G_\K(p;\infty) = |\Phi(\wh p)|$, where $\wh p \in \wh \CC$ is any of the two points above $p\in \CC$.
\item 
The set $\K$ has the $S$--property, namely, 
at each smooth point of $\K$ the opposite normal derivatives of $G_\K$ coincide:
\be
\frac {\pa G_\K(p)}{\pa {\bf n}_+} = \frac {\pa G_\K(p)}{\pa {\bf n}_-}
\label{Sproperty}. 
\ee
\item 
Finally, $\mathfrak Q$ is a Boutroux quadratic differential, namely, all closed contour integrals  of $v=\sqrt {\mathfrak Q}$ on the double cover $\wh \CC$:
\be 
\label{Boutroux}
\oint _\gamma v\in i\R, \ \ \ \forall \gamma \in H_1(\wh C\setminus \{\wh \infty_\pm\}),
\ee
where $H_1$ denotes the homology group of $\wh C$.
\end{enumerate}
\ec

\noindent {\bf Proof.}
For the proof of the first point we could follow Lemma 5.3 of \cite{MFRakh} or Proposition  3.8 of \cite{KuijlaarsSilva} since all the considerations therein are of local nature. We provide also an alternative way of more geometrical nature.

 By definition of Green's function $\K = G_\K^{-1}(\{0\})$ is the zero level-set; first of all  the Green's function is non-trivial only on the  outer region $\Omega$. Different poly-continua that have the same $\Omega$ have thus the same capacity, and we recall our assumption \eqref{ass} so that we can assume that $\K = \CC \setminus \Omega$. Now, $\K$ is must have no interior as it is easy to see that the capacity is monotonic under inclusion, and any disk contained in $\K$ could then be replaced by an arc thus reducing the capacity and preserving the connectivity pattern. Thus the minimal $\K$ should be a poly-continuum without interior points and such that $\Omega$ has full measure. 
 
 Now, since $2\pa_p G_\K(p) = \sqrt{\frak Q}$, we conclude that the square root can be defined consistently across $\Omega$ to a meromorphic differential with a unit  bi-residue at $\infty$. 
Let us understand then $v = \sqrt{\frak Q}$ as a meromorphic differential in $\CC\setminus \K$; in particular any simple pole of $\frak Q$ and any odd-multiplicity zero must then belong to $\K$ because otherwise it would not be possible to define the square-root as a single--valued differential in the neighbourhood of such a  point. It then follows that 
\begin{itemize}
\item For any closed loop $\gamma\in \CC\setminus \K\cup\{\infty\}$ we have $\oint_\gamma v \in i\R$. 
\item We can write $G_\K(p) = \Re \int_{p_0}^p \sqrt{\frak Q}$, where the integral is any path joining $p_0$ to $p$ in the above domain (the integral is defined up to additive imaginary constant by the previous point). 
\end{itemize}
Since $G_\K$ coincides with the above integral, defining a harmonic function $\phi = \Re \int_{p_0}^p \sqrt{\frak Q}$  on $\CC\setminus \K\cup\{\infty\}$, it follows that $\K$ must coincide with the level set of 
$\phi$ from $p_0$, namely, a union of vertical trajectories of $\frak Q$ \cite{StrebelBook}; in particular $\K$ must be piecewise smooth.  This implies the $S$-property since the differential $2\pa_pG_k=  \sqrt{\frak Q}$ extends analytically in the neighbourhood of any smooth point of $\K$. 

Let us define $\wh \CC$ as the double-cover of $\CC$ obtained by glueing two copies along the set $\K$; on $\wh \CC$ the function $\phi$ extends to a function $\Phi(\wh p)$ where $\Phi(\wh p)$ is $-\phi(p) $ for $\wh p$ on the second sheet. This also means that on the second sheet $\sqrt{Q}$ takes the opposite sign. This function $\Phi$ is then harmonic on $\wh \CC$ an odd under the involution. In particular then $2 \pa_{\wh p} \Phi = v$ on both sheets and since $\Phi$ is harmonic, the Boutroux condition \eqref{Boutroux} follows. In conclusion we have $G_\K  = |\Phi|$ by construction.  \QED
We will see that we have also a converse of the Theorem \ref{quadrdiff} and Corollary 
\ref{corboutroux}. 
Namely we can state the following theorem:
\bt
\label{converseBoutroux}
Let $\mathfrak Q$ be a meromorphic quadratic differential  on  $\CC$ with at most simple poles at a finite divisor $E$ and a single double pole at $\infty\in \CC$ with bi-residue $1$. Suppose also that it satisfies the {\it Boutroux condition} \eqref{Boutroux}. 
Then there is a harmonic function $\Phi:\wh \CC \setminus \{\infty_\pm\}\to \R$ such that $2\pa_{\wh p}\Phi =  v$  and  $\Phi(\wh p^\star) =- \Phi(\wh p)$. 
Moreover the function  $G_\K(p):= |\Phi(\wh p)|$ is the Green function of the set, $\K = \pi \le(\Phi^{-1}(\{0\})\ri)$, consisting of critical trajectories of $\mathfrak Q$, and this set is minimizing the capacity within the class with the same connectivity pattern. Furthermore $\K$ has the $S$-property \eqref{Sproperty}. 
\et
\noindent
{\bf Proof.}
Recall that ${}^\star:\wh \CC\to \wh \CC$ is the holomorphic involution exchanging the sheets and that 
the differential $v = \sqrt{\mathfrak Q}$ on $\wh \CC$, by definition, is odd under this involution. 
Consider the integral 
$
\int_{\wh p^*}^{\wh p} v$: the contour of integration is defined up to the homology of $\wh \CC \setminus \{\wh \infty_\pm\}$ and hence the integral is defined up to the periods of $v$, which are purely imaginary by assumption. 
Thus 
\be
\Phi(\wh p):= \Re\le(\int_{\wh p^\star}^{\wh p} v\ri) 
\ee
is a well--defined harmonic function on $\wh \CC \setminus \{\wh \infty_\pm\}$ which is manifestly odd under the involution, $\Phi(\wh p^\star)=-\Phi(\wh p)$. This function behaves like $\mp \ln |z_\infty(\pi(p))| + \mathcal O(1)$ as $\wh p\to\wh \infty_\pm$. Let $\lambda\in \R$ and $E(\lambda):= \Phi^{-1}(\{\lambda\})$ denote the $\lambda$--level set. Since   $\Phi$ is odd under the involution we must have   $E(\lambda)) = E(-\lambda)^\star$,   and hence it makes sense to consider $\pi(E(\lambda))$ as a function of $|\lambda|$. For $\lambda=0$ we then have $E(0)=E(0)^\star$ and hence this  is a well-defined set of $\CC$, $\K\subset \CC$. 
The function $G_\K(p):= |\Phi(\wh p)|$, with $p=\pi(\wh p)$, is therefore a positive continuous function in $\CC\setminus \{\infty\}$, harmonic in $\CC\setminus \{\infty\}\cup \K$ and behaving like $\ln |z_\infty|$ near $\infty$. Thus it is, by definition, the polar Green's function of the complement of $\K$. 

Since $E(0)$ is the level-set of a harmonic function $\Phi$, it consists of a finite union of analytic arcs: the gradient lines of $\Phi$ are orthogonal to $E(0)$ almost everywhere (where $E(0)$ is smooth), and thus the normal derivatives of $G_\K$ on the two sides of $\K=\pi(E(0))$ are equal  to each other, $\pa_{{\bf n}_+}G_\K =\pa_{{\bf n}_-} G_\K$ since each of them is the absolute value of the gradient of $\Phi$. In other words, the set $\K$ in $\CC$ has the $S$--property. 

The proof of the minimization of the capacity amongst the class, $\mathbb K$, of all polycontinua in the same connectivity class is contained in Theorem \ref{thmJIP} and Corollary \ref{corS}.  This is a consequence of the fact that any such set has the Jenkins interception property relative to $\K$. 
\QED

\section{The generalized Jenkins Interception Property}
\label{secJIP}
In this section we formulate and prove a criterion that can be used to establish an inequality between the capacities of two poly-continua $\frak F,  \K$. The criterion is a generalization of the result of \cite{JenkinsArt}, which was proved for poly-continua of the Riemann sphere that  satisfy what is now known as the $S$-property \cite{Stahl, Stahl2,Stahl3, MFRakh}.

Let $\frak F\subset \CC$ be a poly-continuum not containing $\infty$: we need to assume that the boundary is piecewise smooth. Let $\Omega$ be the outer region of $\frak F$.  Let  $G_{\frak F}(p;\infty)$ denote  the corresponding Green's function.  
All but finitely many gradient lines of $G_{\frak F}$ connect the outer boundary $\pa \Omega$  to $\infty$. Following Jenkins, we denote with $\mathfrak L$ the set of gradient lines of $G_{\frak F}$ (the ``orthogonal trajectories'' to the level-set foliation) that extend from $\pa\Omega$ to $\infty$.   This is also a {\it measured foliation}, where the measure is determined by the variation of the conjugate of the Green's function. To be more specific, each $\gamma \in \mathfrak L$ originates from $p\in \pa\Omega$ and the density of the measure (with respect to the arclength) is $\pa G_{\frak F}/\pa {\bf n}$, with ${\bf n}$ the normal pointing to the interior of $\Omega$ along the direction of $\gamma$. The total mass of this measure is $2\pi$.  Note that if $\frak F$ contains arcs without interior, there are two trajectories in $\mathfrak L$ from each smooth point of such arc. 

\bd
Let $p\in \pa \Omega$ be a smooth point of the boundary;  if there are two trajectories originating at $p$ in opposite directions ${\bf n}_{1}=-{\bf n}_2$
 we will say that the {\bf dominant} trajectory, $\mathcal L^\perp_d(p)$, is the one for which $\pa_{{\bf n}_j} G_{\frak F}(p)$ is larger, while the opposite trajectory $\mathcal L^\perp_r(p)$ will be called {\bf recessive}.
 If both normal derivatives are equal (S-property), then both trajectories are declared dominant. 
  If there is only one trajectory\footnote{This happens if the boundary locally bounds a region of $\frak F$ with non-empty interior.} it will be the dominant by default.  We denote $\mathfrak L_d$ the collection of the dominant trajectories. 
\ed

With this notion we can formulate the generalized Jenkins interception property as follows.
\bd
 We say that a poly-continuum $\K$ has the {\bf Jenkins Interception Property} (JIP) relative to $\frak F$ if $ \K$ intercepts all dominant trajectories:
\be
\forall \gamma\in \mathfrak L_d \ \ \gamma\cap  \K \neq \emptyset.
\ee
\ed
We warn that the interception may occur at the point $p\in \pa \Omega$. 
We can now state the comparison criterion. 
\bt
\label{thmJIP} 
Suppose $\frak F$ is as described in this section and $ \K$ is another poly-continuum satisfying the Jenkins Interception Property relative to $\frak F$. Then 
\be
{\rm Cap}_{z_\infty}(\frak F) \leq {\rm Cap}_{z_\infty}( \K)
\ee
The equality is realized if and only if $ \K = \frak F$
\et
We recall our running assumption \eqref{ass}, namely, that all poly-continua we mention are tacitly assumed to be polynomially convex. Without this assumption the equality mentioned in the theorem could fail in trivial ways. 

\noindent{\bf Proof.}
The proof follows almost verbatim that of Jenkins in \cite{JenkinsArt}  but with a bit more details for the sake of the reader, see also \cite{BertoTovbisNLSCritical}. 
We write $G(p):= G_{\frak F}(p)$ for brevity unless the specification of the set is important for clarity. We denote $w(p):= {\rm e}^{-g_\frak F}$ as in Section  \ref{secunif} and defined on the surface minus the cuts as explained there. 
We denote by $W(p) = {\rm e}^{-g_\K(p)}$ instead.  The inverse map to $w(p)$ will be denoted $\p(s,\theta)$ with $w = {\rm e}^{-s+i\theta}$. This inverse map is defined on the strip $s\in [0,\infty)$, $\theta \in [0,2\pi)$ (with periodic  identification) up to a finite number of vertical slits of finite length in the metric $\sqrt {\d s^2  + \d \theta^2} = |\pa_p G|$.
  
For $p$ near $\infty$ we have the relation 
\be
\label{Ww}
W(p)  = \CC\, w(p)\bigg(1 + \mathcal O(w)\bigg), \ \ \ \CC:= {\rm e}^{\I(\K) - \I(\frak F)}.
\ee
Note that since both $\I$'s in the exponent are defined up to addition of the same constant (depending on the normalization of the bipolar Green function), the difference is well defined. We need to show that the difference is negative.
Let $\Delta(\l;\frak F):= G^{-1}\big([0,\l]\big)$. 
For the set $\K$ we choose $\wt \l$ to be  the largest value for which $G_\K^{-1}([0, \wt\l])=:\Delta(\wt\l;\K) \subset \Delta(\l;\frak F)$. 
From \eqref{Ww} we deduce that for $\l$ large enough the values $\wt \l ,\l$ are related as follows:
\be
\wt \l = \l + \delta + \mathcal O(\l^{-1}) ,\ \ \ \delta:= \I(\K)-\I(\frak F). \label{wtll}
\ee
Let $R\subset [0,\l]\times [0,2\pi)$ be the image of the region $\Delta (\wt \l\;\K)\subset \Delta(\l;\frak F)$ under the map $(\ln|w|,\arg(w))$. 
Let us define the conformal factor in $\Delta(\l;\frak F)$ of the metric $\frac{|\pa_p G_\K|}{|\pa_p G|} = \rho \sqrt{\d s^2+\d \theta^2}$: 
\be
\rho(s,\theta):=\le| \frac {\d g_\K}{\d g_{\frak F}}\ri| \mathcal X_{R},
\ee
where $\mathcal X_S$ denotes the indicator function of a set $S$. This conformal metric is the pull-back to $[0,\l] \times [0,2\pi)$ of the metric $|\pa_p G_\K|$ induced by the Green's function of $\K$ on $\Delta(\wt \l;\K)$.  
We give an estimate from below of the lengths  of the gradient lines of $G$  in the $ \rho$--metric: these are the curves $\p(s,\theta)$ for $\theta\in [0,2\pi)$ fixed and $s$ ranging in $[0,\l]$ and we need to estimate 
\be
\L_{ \rho}(\theta):=  \int_0^{\l} \rho(s,\theta) \d s .
\ee

Let $\Omega$ be the  outer region of $\frak F$  and
let  $\frak F^o$ denote the smooth part of the outer  boundary  $\pa\Omega$. Each $p\in \frak F^o$ is accessible from at least one side from $\Omega$ and $w(p) ={\rm e}^{i\theta}$. If $p$ is accessible also from the other side, we denote it by $p^\star$ and correspondingly $w(p^\star) = {\rm e}^{i\theta^\star}$. In this latter case the dominant point is the one for which the limit of the gradient of $G$ along the gradient line ending at that point is greater amongst the two values. If the limits are equal, both $p, p^\star$ are declared dominant.  
 
 Let us define the {\it dominant subset} by 
 \be
 U^+:= \{\theta\in [0,2\pi): \ \text {the trajectory $\p( [0,\lambda],\theta)\subset \CC\setminus \frak F$ is dominant}\}.
 \ee
 Similarly the {\it recessive subset} is defined by 
 \be
 U^-:= \{\theta\in [0,2\pi): \ \text {the trajectory $\p( [0,\lambda],\theta)\subset \CC\setminus \frak F$ is recessive}\}.
 \ee
 Recall that the points  of $\frak F^o$ accessible from only one side are dominant by definition. 

We claim that 
for $\lambda$ sufficiently large we have 
 \bea
\label{estmain-adj}
\mathcal L_\rho(\theta) = &\int_{0}^\l  \rho(s,\theta)\d s \geq
  \wt \l -  {G_\K (p)} ~~\text{when }~~\theta \in U^-
\cr
\mathcal L_\rho(\theta)=&\int_{0}^\l  \rho(s,\theta)\d s \geq   \wt \l +  {G_\K (p)} ~~\text{when }~~\theta\in U^+
\eea
Here $p =\p(0,\theta)\in \frak F^o$. \\
\noindent {\bf Proof of claim.}
In the first case the integral is the total variation of $ g_\K$ along the gradient line of $G$ that starts at $p$ and ends at $\p(\lambda, \theta)$. This is greater than the variation of $G_\K = \Re(g_\K)$. Since the geodesic exits $\Delta(\wt \l;\K)$ by construction of $\wt \l$ (see lines above \eqref{wtll}), this variation is at least $\wt \l - G_{\K}(p)$. 
Consider now the second case: the integral is still the total variation of $ {g_\K}$, which is greater than the total variation of $G_\K$ exactly as before. However we have assumed that $\K$ intercepts all dominant trajectories and hence $G_\K$  moves from $G_K(p)>0$ down to zero (when the trajectory meets $\K$) and then from $0$  up to $G_\K(\p(\l,\theta))>\wt \l$. Since we chopped the metric to $\Delta(\wt \l; \K)$ the maximum value will be $\wt \l$.   \hfill $\blacktriangle$

We recall that, by our definition, for each recessive trajectory there is a dominant one (but not necessarily the converse), and hence $(U^-)^\star \subseteq U^+$ under the map $\theta\mapsto \theta^\star$. 
The pull-back of the Lebesgue measure $\d\theta$ from $U^-$ to $(U^-)^\star$ is 
\be
\phi(\theta)\d\theta = \frac {\pa_{\bf n_-}G}{\pa_{\bf n_+} G} \leq \d\theta, \label{dm2}
\ee
where ${\bf n}_+$ is the normal pointing in the dominant and ${\bf n}_{-}$ the one pointing in the recessive direction.

Now we need another metric defined on $\Delta(\l;\frak F)$; this is a partial average of $\rho$ as follows.
Now define 
\be
\wh \rho(s,\theta):= \frac {\rho(\theta,s)\mathcal X_{U^-}  + \rho(\theta, s)\mathcal X_{(U^-)^\star} }2
+
 \rho(\theta,s) \mathcal X_{U^+\setminus (U^-)^\star}
\ee
The two estimates \eqref{estmain-adj} imply that the length of $\p([0,\l],\theta)$ in this metric $\wh \rho$ satisfies
\be
\label{612}
\mathcal L_{\wh \rho}(\theta) = \int_0^\l \wh \rho(s,\theta)\d s\geq \wt \l, 
\ee
for almost all $\theta$'s.  This means that $\wh \rho/\wt \l$ is an admissible metric for the module problem for $\Delta(\l,\frak F)$ \cite{StrebelBook}, Ch. II.
In particular this means that 
\be
\scr A_{\frac{\wh \rho}{\wt \l}} ([0,\l]\times [0,2\pi) ) \geq \frac {2\pi}{\l}, \label{est11}
\ee
since the latter is the module in the extremal metric $\frac{\sqrt{\d s^2+\d\theta^2}}\l$.
On the other hand, with $\phi$ defined in \eqref{dm2}, we have 
\bea
\scr A_{\wh \rho} ([0,\l]\times [0,2\pi) ) \leq \frac 1 2 \int_{U^-} \!\!\! \int_0^\l \rho^2(s,\theta) \d s\d\theta +
 \frac 1 2 \int_{(U^-)^\star} \!\! \int_0^\l \rho^2(s,\theta) \d s\d\theta\cr 
 +    \ \int_{U^+\setminus (U^-)^\star} \!\!\int_0^\l \rho^2(s,\theta) \d s\d\theta=\cr
\leq 
\ \int_{(U^-)^\star} \!\! \int_0^\l \frac {1+\phi^2(\theta)}2 \rho^2(s,\theta) \d s\d\theta 
 +    \ \int_{U^+\setminus (U^-)^\star} \!\!\int_0^\l \rho^2(s,\theta\d s\d\theta )\leq 
 \cr 
 \leq   \int_{U^+} \!\!\int_0^\l \rho^2(s,\theta) \d s\d\theta \leq \int_0^{2\pi}  \!\!\int_0^\l \rho^2(s,\theta) \d s\d\theta = \scr A_\rho([0,\l]\times [0,2\pi) )
 \label{upper}
\eea
This latter integral, due to the presence of the characteristic function of $\Delta(\wt \l;\K)\subset\Delta(\l;\frak F)$, computes precisely the module of $\Delta(\wt \l;\K)$, which is,  the area of the annulus $|\zeta|>{\rm e}^{-\wt \l}$ in the metric $\frac {|\d \zeta|}{\wt \l|\zeta|}$ in the $\zeta = W_\K$ plane.  Namely
\be
\scr A_\rho([0,\l]\times [0,2\pi) )= \frac {2\pi}{\wt \l}.
\ee

Putting together the estimates we  have shown 
\be
\frac {2\pi}{\l} \mathop{\leq}^{\eqref{est11} } \scr A_{\wh \rho} (\Delta(\l;\frak F))  \mathop{\leq}^{\eqref{upper}} \frac {2\pi }{\wt \l}  \label{616}
\ee
from which it follows 
\be
\wt \l = \l+\I(\K)-  \I(\frak F) + \mathcal O(\l^{-1}) \leq \l \ \ \Leftrightarrow \ \I(\K)\leq \I(\frak F) \ \Leftrightarrow \ {\rm Cap}_{z_\infty}(\K)\geq {\rm Cap}_{z_\infty} (\frak F). 
\label{617}
\ee
This establishes the desired inequality. 

It remains to explain how the equality can be achieved. The argument here follows \cite{JenkinsArt}, although Jenkins  is extremely sparse in his explanation and for us to be  more transparent, we need to unpack a bit the ``length-area method''. 
If the capacities are equal, we have from \eqref{617} that $\wt \l=\l + \mathcal O(\l^{-1}) $ and hence,
\bea
 \frac 1{\wt \l^2}  \iint \le( { \rho^2}-1 \ri) \d s \d\theta-   \frac {2\pi}{\l}+ \frac {2\pi \l}{{\wt \l}^2} 
=
 \iint \le(\frac{ \rho^2}{\wt \l^2} - \frac 1{\l^2} \ri)\d s \d\theta = \frac {2\pi}{\wt \l} - \frac {2\pi}\l,
 \eea
 namely 
 \bea
 \int_0^\l\int_0^{2\pi} \le( { \rho^2}-1 \ri) \d s \d\theta=\mathcal O(\l^{-1}).
 \eea
 Now,
 \bea
 \mathcal O(\l^{-1}) =  \iint \le( { \rho^2}-1 \ri) \d s \d\theta= \iint \le( { \rho}-1 \ri)^2 \d s \d\theta +2  \iint \le( { \rho}-1 \ri) \d s \d\theta .
 \label{qq2}
 \eea
Now we have from \eqref{estmain-adj}
\bea
\int_0^{2\pi} \L_{{\rho} }(\theta) \d \theta = 2\pi \wt \l  +\int_{U^+} G_\K(p(\theta)  \d \theta - \int_{U^-} G_\K(p(\theta) \d \theta =2\pi \wt \l +\int_{U^+} G_\K(p(\theta)(1-\phi(\theta))   \d \theta \geq 2\pi \wt \l, 
\eea
which can be rewritten as 
\bea
\label{qq1}
\iint(\rho-1)\d s \d \theta \geq  2\pi \wt \l - 2\pi \l = \mathcal O(\l^{-1}).
\eea
Using \eqref{qq1} into \eqref{qq2}  we have 
\be
 \int_0^\l\int_0^{2\pi} \le( { \rho}-1 \ri)^2 \d s \d\theta = \mathcal O(\l^{-1}), \label{mw1}
\ee
 so the only possiblity for \eqref{mw1} to hold is that $\rho\equiv 1$. 
Recall now that $\rho = \le|\frac {g'_\K}{g'}\ri|$ and that both $g_\K', g'$ are locally analytic in a punctured disk centered at  $\infty$. Since an analytic function with constant modulus must be constant, and since both  $g_\K , g$ behave like $-\frac {\d z_\infty}{z_\infty}$ we conclude that $g'_\K \equiv  g'$  where they are both defined. Namely $g_\K, g$ differ by a constant in the common domain. But since both are supposed to be green's functions, we conclude that $\K = \frak F$. \QED

\bc
\label{corS}
Suppose that $\frak F$ is the set associated to a Boutroux quadratic differential as in Theorem \ref{converseBoutroux}.  Then all but finitely many orthogonal trajectories are dominant  and $\frak F$ is the unique minimizer in the class of polycontinua with the same connectivity pattern. 
\ec
\noindent{\bf Proof.}
Since $\frak F$ is a finite union of analytic arcs at all smooth points there are two trajectories and  we know already from Theorem \ref{converseBoutroux} that $\frak F$ has the $S$-property \eqref{Sproperty}. Thus all but finitelly many trajectories are dominant. Let $\K$ be a polycontinuum in the same connectivity class. We need to show that it satisfies  the Jenkins Interception Property relative to $\frak F$ and hence has greater capacity. Let $\scr P$ denote  the connectivity pattern determined by $\frak F$ in the following way; for every $e,\wt e\in E$ that belong to the same connected component of $\frak F$ there is a minimal union of analytic  sub-arcs of $\frak F$ connecting the two, and this determines the element in $\pi_1(\CC,E)$ used to define ${\scr P}$. We denote this path by $\gamma_{e,\wt e}$. . Let $\K\in \mathbb K_{\scr P}$ be any other polycontinuum. By definition of the class $\mathbb K_{\scr P}$, for every $\epsilon$ fattening, $\K^{(\epsilon)}$ of $\K$  and for each $e,\wt e\in E\subset \K$ there is thus a path $\gamma_{e,\wt e}^{(\epsilon)}$ that is homotopic to $\gamma_{e,\wt e}$. These two paths bound a simply connected region that is foliated by the gradient lines of the Green's function $G_{\frak F}$. Thus for all the smooth point $p$  of $\gamma_{e,\wt e}$ at least one of the two gradient lines must intersect $\eta_{e,\wt e}^{(\epsilon)}$. 
Since $\epsilon>0$ is arbitrary, a subsequence of these  points of intersection must converge, as $\epsilon \to 0$, to a point of $\K$ and hence $\K$ has the Jenkins interception property. 
Then by Theorem \ref{thmJIP} $\frak F$ must be a minimizer within $\mathbb K_{\scr P}$: the same theorem also shows that this is the unique minimizer. 
\QED
To tie the loose ends we still need to prove the uniqueness of the minimizer of Theorem \ref{thmP}.
\bp
\label{propuniqueness}
For any admissible connectivity pattern $\scr P$ the minimum of the capacity over the class $\mathbb K_{\scr P}$ is achieved at a {\rm unique}  set $\K_{\scr P}$.
\ep
\noindent{\bf Proof.}
By Corollary \ref{corboutroux} the minimizer $\K_{\scr P}$ satisfies the conditions of Corollary \ref{corS}  and hence it is the unique minimizer in the class.
\QED
\section{Jenkins-Strebel differentials and minimal capacity}
\label{CommentJS}
The purpose of this section is to comment on,  and elucidate the relationship between the notion of Jenkins-Strebel differential and problems of minimal capacity, which  was already pointed out in \cite{BertoSAPM} for the case of the Riemann sphere. 
We refer to the review article \cite{Arbarello} for many of the relevant notions. 
A Jenkins-Strebel  (JS) differential $\frak  Q$ is a quadratic differential on a Riemann surface $\CC$ with at most simple poles at a set $E=\{e_1,\dots, e_n\}$ and with double poles at another set $P = \{\infty_1,\dots, \infty_k\}$ satisfying the following properties (see Def. 31 in \cite{Arbarello}):
\begin{enumerate}
\item The bi-residues at each of the double poles are positive, namely, in local coordinate $z$ near $p_j$
\be
\frak Q =  \frac {a_j^2}{4\pi^2} \frac {\d z^2}{z^2} (1 + \mathcal O(z)),\ \ a_j>0.
\ee
\item All non-critical vertical trajectories are closed loops.
\end{enumerate}
For ease of comparison with the current setting, we have changed the overall sign of the quadratic differential in the definition of \cite{Arbarello}, so that the horizontal trajectories therein are now vertical.
Then the Jenkins-Strebel theorem is as follows. 
\bt[Theorem 32 in \cite{Arbarello}]
\label{JSdiff}
For any sets $E, P$ as above, and for any choice of parameters $a_j>0$ there is a unique JS differential such that the critical graph  (the union of the critical trajectories) is connected and contains $E$, and the complement of the critical graph is  a disjoint union of conformal disks $\DD_j$, each containing exactly one of the points $p_j$. 
\et
Jenkins-Strebel differentials were used to construct a cellular decomposition of the moduli space of Riemann surfaces of genus $g$ with $k$ marked points (these are the set $P$) in unpublished works of Mumford and Thurston, and by Harer \cite{Harer1,Harer2}.  The critical graph is called also {\it ribbon graph} and features prominently in the proof by Kontsevich of the Witten conjecture \cite{Kontsevich, Witten}
. 
It is an instructive exercise which we leave to the reader  to show that the Jenkins-Strebel differentials of Theorem \ref{JSdiff} are Boutroux differentials and, as a partial converse, a quadratic differentials with at most double poles, satisfying the Boutroux condition and with simple zeros has the properties in Theorem \ref{JSdiff}. Note, however, that  just the property of having closed vertical non-critical trajectories is weaker than the Boutroux condition.

The purpose of this section is rather to elucidate the relationship with our Chebotarev problem. 
Suppose that the connectivity pattern $\scr P$ for the set of anchors is such that the complement of the minimizing set $\K$ is simply connected. Note that, in genus $0$ this case is the classical Chebotarev problem of finding a continuum containing the anchors, but in higher genus, in general,  the mere requirement that the minimizing  continuum contains the anchors does not imply that the complement is simply connected. 

Suppose that $\frak Q$ is the quadratic differential in our Theorem \ref{quadrdiff} for the minimizer in the class of the above connectivity pattern; then $\frak Q $ is the Jenkins-Strebel differential of Theorem \ref{JSdiff} with $P=\{\infty\}$ and $a_1 = 2\pi$. In other terms $\K$ is the ribbon graph of Harer. 
 
 Viceversa, the Harer ribbon graphs in the case of one double pole are sets of minimal capacity in our sense. Of course for the case of several double poles  what is being minimized is the sum of the ``reduced modules'', see \cite{StrebelBook}, where the reduced module of a disk is equivalent to the notion of capacity here. 
 
 We do not delve further into this connection but we wanted to point it out because we have noticed that the potential theory community and the moduli space community  by and large seem to have not realized the overlap between these two problems. 

\section{Conclusion and other similar problems}
\label{conclusion}

We want to comment on the application of this result to the study of Pad\'e\ approximation and orthogonality. We keep this sketch of example ``simple'' for the sake of discussion, but the reader with experience in the asymptotics of orthogonal polynomial may easily come up with formulations of more general problems including  those with scaling weights.

Let $E = \{e_1,\dots, e_N\}$ and choose  paths connecting these points  together in some connectivity pattern. Denote by $\Sigma = \bigcup_{j=1}^L \sigma_j $ the union of these paths.  
Fix also  a point $\infty\not\in \Sigma$ and a differential $\eta$ defined and holomorphic on the universal cover of $\CC\setminus E\cup\{\infty\}$ with the following local behaviour near the points
\bea
 \eta &= z_j^{\alpha_j} (C_j+\mathcal O(z)) \d z_j, \ \ \Re \alpha_j \geq -1,\ \ 
\nn\\
 \eta &= z_\infty^{-\sum \alpha_j} (C_\infty + \mathcal O(z_\infty)) \d z_\infty, \nn
\eea
where $z_j, z_\infty$ are local coordinates near $e_j$ and $\infty$. 
Following \cite{BertoPade} we fix a generic divisor $\scr D$ of degree $g$ not intersecting $\Sigma$ or $\infty$.
It follows from the Riemann Roch theorem that the linear space of meromorphic functions with pole at $\scr D$ and at most order $n$ at $\infty$ are a space of dimension $n+1$, much like the polynomials of degree not exceeding $n$.
Then we can pose the problem of finding a sequence $(\pi_n)$ of  meromorphic  functions with at most simple poles at $\scr D$ and pole of order $n$ at $\infty$ such that 
\be
\label{pipieta}
\sum_{j=1}^{L} \varkappa_j \int_{\sigma_j} \pi_n \pi_m \eta = \delta_{nm},
\ee
where $\varkappa_j$ are arbitrarily chosen complex coordinates.
This  is  a non-hermitian orthogonality of sorts; the setting mimicks the situation of Jacobi polynomials that satisfy 
\be
\label{jacpol}
\int_{-1}^1 P_n(x) P_m(x) (x+1)^\alpha (1-x)^\beta \d x = \delta_{nm}.
\ee
Here $\CC = \mathbb P^1$, $E = \{1,-1\}$, $\infty$ is the point at infinity and $ \eta =(x+1)^\alpha (1-x)^\beta \d x $. In this simple example, the contour of integration in \eqref{jacpol} is any arc connecting the two anchors, but the zeros of the orthogonal polynomials irrespectlvely of the choice of contour will accumulate on the straight segment, which is the minimal capacity set in this class. 

Following \cite{BertoPade} and the template of \cite{BertoEllOPs} one can setup a Riemann-Hilbert problem and study the large $n$ asymptotics for the orthogonal meromorphic functions $\pi_n$. In \cite{BertoPade} it is explained how the  $\pi_n$'s appear as denominators of the one-point Pad\'e\ approximation problem of a suitable Stieltjes transform of the measure $\eta$ on $\Sigma$ . 

The first step in the nonlinear steepest descent analysis requires the construction of a $g$--function (in the lore of the literature) and then  follows  steps similar to \cite{DKMVZ}. The $g$ function for the described problem would be exactly the function $g_\K = G_\K + i H_\K$, related to the minimal capacity set $\K$ in the connectivity class determined by $\Sigma$. 
This example should suffice to explain our interest in the problem discussed in this paper. 
\vspace{-12pt}
\paragraph{Generalizations.}
\begin{itemize}
\item The Green's function $G_\K$ takes zero value on $\K$. However one could consider a problem where on different components the function takes different constant values and then study the minimal ``capacity'' of these problems, with the capacity being defined by the expansion at infinity as per Definition \ref{defGreenK}. 

Let us outline very briefly how this problem would be relevant in a similar context. 
Suppose now we are studying the large $n$ degree of orthogonal functions for the pairing \eqref{pipieta} but  suppose now that the constants $\varkappa_j$ are allowed to scale at an exponential rate with $n$, namely, we consider a problem with varying weight of rather simple form. The analogy for ordinary orthogonal polynomials would be the study of the following  orthogonal polynomials, for example 
\be
\int_0^1 P_n P_m \d x + {\rm e}^{-N}\int_3^4  P_n P_m \d x,
\ee
where we study the behaviour of $P_n$ as $n\to \infty$ and $\frac N n\to c\neq 0$.
The $g$ function for this problem would then be constructed out of the obstacle  problem where $G_\K$ takes different  values greater or equal to  $\frac 1 n \ln \varkappa_j$ on the separate components. Note that the relevant problem to be solved is really the obstacle problem rather than the Dirichlet problem; however the two are equivalent in the setting of this paper.

\item Of course the full generalization would be to consider an arbitrary scaling weight, for example of the form ${\rm e}^{-n V(z)} \eta$ with $V(z)$ meromorphic.  Then, instead of a pure capacity problem we would be lead to the study of a potential problem with external field. This type of problems become rapidly difficult to handle but in the case of ordinary orthogonal polynomials the energy max/minimization problem was successfully tackled at least for certain external fields, \cite{MFRakh, KuijlaarsSilva} and for $\CC$ of genus one some progress was in \cite{BertoGrootKuijlaars}.

\item Our motivation stems from the one-point Pad\'e\ approximation problem. If we had several points of simultaneous approximation, $\infty_1,\dots, \infty_k$ with the corresponding orders of approximation, $n_1,\dots, n_k$ growing at the same rate,  we would have to study a similar problem but where the Green's function has several logarithmic poles, of different relative positive strengths adding up to  $1$.  This would lead to a generalized Chebotarev problem in the presence of several external charges and make full contact with the theory of Jenkins-Strebel differentials mentioned in Section \ref{CommentJS}.
\end{itemize}
We hope to return to these problems in the near future.

%
%
%
%
%
%
%
\appendix
\section{Explicit formulas for the Cauchy kernel $C^{(2,-1)}$}
\label{expCauchy}
We recall basic facts and notation.
\label{theta}
We choose a Torelli marking $\{a_1,\dots, a_g, b_1,\dots, b_g\}$ for $\mathcal C$ in terms of which we construct the classical Riemann Theta functions. We refer to \cite{Fay}, Ch. 1-2 for a review of these classical notions.  The corresponding normalized Abelian differentials, period matrix and Abel map will be denoted by $\omega_j, \bs \tau, \mathfrak A$, respectively:
\be
\oint_{a_j} \omega_k = \delta _{jk},\  \ \bs \tau_{ij}= \oint_{b_j} \omega_i = \oint_{b_i} \omega_j,\ \ \ \mathfrak A_{p_0}(p)= \int_{p_0}^p \vec \omega: \mathcal C\times \mathcal C \to \mathbb J(\mathcal C).
\ee
 The base-point of the Abel map will not be indicated except when necessary. The Abel map is naturally extended to a map on divisors.
With $\Delta$ we  denote  characteristics $\Delta = [\bs\alpha, \bs \beta]\in \R^g\times \R^g$ corresponding to the point  ${\mathbf e}=\bs \beta - \bs \tau \bs \alpha\in \mathbb J(\mathcal C)$ in the Jacobian.
The Theta function with characteristic $\Delta$ is the function defined on $\C^g$ by 
\be
\Theta_{\Delta}({\bs z};\bs \tau) :=
 \sum_{{\bs n} \in \Z^g} 
  \exp \le[i\pi 
 \le({\bs n} +\frac 1 2  \bs \beta\ri)^t
 \cdot \bs \tau \cdot\le({\bs n}  +\frac 1 2\bs \beta\ri) + 2i\pi \le({\bs n}  +\frac 1 2\bs \beta\ri)\cdot \le({\bs z} + \frac 1 2 \bs  \alpha\ri)
 \ri], \ \ {\bs z}\in \C^g.
\ee
It has the periodicity  properties:
\be
\label{periodprop}
\Theta_\Delta(\bs z + \bs m + \bs \tau \bs k;\bs \tau) = {\rm e}^{2i\pi \le(
\frac {\bs \alpha \bs m - \bs \beta \bs k}2  - \bs k \bs z - \frac 1 2 \bs k\bs \tau \bs k
\ri)} \Theta_\Delta(\bs z;\bs \tau).
\ee
If the characteristics $\Delta$ is made of integers, we say that $\Delta$ is {\it odd} if $\bs \alpha \cdot \bs \beta\in 2\Z+1$: in this case, fittingly, $\Theta_\Delta(-{\bs z};\bs \tau) = -\Theta_\Delta({\bs z};\bs \tau)$.

It is known that there always are {\it nonsingular} odd characteristics (Corollary 4.21 in \cite{Fay}), namely, one for which the gradient $\nabla \Theta_\Delta(\bs 0)$ is not the zero vector. For any such odd nonsingular characteristic, the following holomorphic differential $\omega_\Delta$ has $g-1$ only double zeros at a divisor $2\scr D_\Delta$, with $\scr D_\Delta$ a divisor of degree $g-1$:
\be
\label{defomegadelta}
\omega_\Delta(p) = \sum_{j=1}^g \pa_{z_j} \Theta_\Delta\bigg|_{\bs0} \omega_j(p). 
\ee
Fix a nonsingular odd characteristic $\Delta$: then \cite {Fay} the multi-valued function 
\be
\Theta_\Delta\le(\int_q^p \bs \omega\ri), \ \ \bs \omega:=[\omega_1,\dots, \omega_g],
\ee
has, as a function of $p$, $g-1$ zeros exactly at $\scr D_\Delta$, independent of $q$. Since $\Theta_\Delta$ is an odd function, the same holds with respect to $q$. 
Furthermore, for any (generic) ${\bs f}\in \C^g$ the multi-valued function $\Theta(\mathfrak A(p)-{\bs f})$ has $g$ zeros (counted with multiplicity).  Finally we recall that (Jacobi Inversion Theorem) for any (generic) divisor of degree $g$ there is a corresponding ${\bs f}\in \C^g$ such that the zeros of the above Theta function coincide with the desired divisor. With these preparations we can provide the expression of the Cauchy kernel $C^{(2,-1)}$ in Section \ref{secminimal}.
\bp
Let $\scr B_{1,2,3}$ be three generic divisors of degree $g$ and ${\bs f}_{1,2,3}\in \C^g$ the corresponding vectors described above. Then 
\bea
C^{(2,-1)}(x,q) = \frac {\Theta_\Delta^{N-3}(\int_\infty^x\bs \omega)}{\Theta_\Delta^{N-3}(\int_\infty^q\bs \omega)}
\prod_{j=1}^N \frac  {\Theta_\Delta(\int_{e_j}^q\bs \omega)}{\Theta_\Delta(\int_{e_j}^x\bs \omega)}
\prod_{j=1}^3 \frac {\Theta(\frak A(x)-{\bs f}_j)}{\Theta(\frak A(q)-{\bs f}_j)}
\frac {\Theta(\int_q^x - {\bf s})}{\Theta_\Delta(\int_{q}^x \bs \omega)} \frac {\omega_\Delta^2(x)}{\Theta(\bs s)\omega_\Delta(q)},
\eea
where $\bs s =\sum_{j=1}^N \frak A(e_j) -{\bs f}_1 - {\bs f}_2- {\bs f}_3 - (N-3)\frak A(\infty)$.
\ep
The genericity of ${\bs f}_{1,2,3}$ guarantees that the constant in the  denominator does not vanish, $\Theta({\bs s})\neq 0$. 
It is an exercise using the periodicity properties \eqref{periodprop} and the other properties mentioned above to verify that the proposed expression satisfies the requirements set in Section \ref{secminimal}. 
In particular, one needs to verify that there are no additional poles beyond the prescribed ones, and that the expression is single--valued on $\CC$, using the periodicity properties \eqref{periodprop}.

\section{Proof of Hausdorff continuity of ${\rm Cap}_{z_\infty}$}
\label{AppProofcont}
The proof proceeds from the following estimate of the H\"older continuity up to the boundary of the Green's function, see \cite{Chirka} Section 2.8; these estimates are based on prior estimates valid in the plane and obtained in an unpublished work  Perevoznikova and Rakhmanov \cite{RakhPer}. See also \cite{Rakhman} and the appendix in \cite{BertoTovbisNLSCritical}.

Let $\mathbb K( \delta, \rho)$ be a class of  poly-continua with the properties that 
\begin{itemize}
\item  there is $\delta>0$ such that for each $\K\in \mathbb K$ the diameter of each connected component of $\K$ is at least $\delta$;
\item there is $\rho>0$ such that $d_\sigma(\K,\infty)\geq \rho$. 
\end{itemize}
Under these hypotheses there is a constant $C = C(\delta, \rho)$ such that 
\be
\label{Greenest}
G_\K(p;\infty) \leq C \sqrt{d_\sigma(p,\K)}. 
\ee
Let us denote by $\K^\epsilon$ the epsilon-fattening of a set $\K$. Suppose that $d_H(\K, \wt \K)\leq \epsilon$. Then $\wt \K\subset \K^\epsilon $ and $\K\subset \wt \K^\epsilon$. But then
\bea
G_{\K}(p;\infty) -\max_{q\in \wt \K} G_{\K(q;\infty)} \leq  G_{ \wt \K}(p;\infty) \ \ \Rightarrow\ \ 
G_{\K}(p;\infty) - G_{ \wt \K}(p;\infty)\leq \max_{q\in \wt \K} G_{\K}(q;\infty) \leq C\sqrt{\epsilon}
\eea
where the last inequality follows from \eqref{Greenest} and the fact that $\wt\K\subset \K^\epsilon$; swapping the roles of $\K, \wt \K$ we obtain then 
\be
\big|G_{\K}(p;\infty) - G_{ \wt \K}(p;\infty)\big|\leq C\sqrt{\epsilon}.
\ee
Taking the limsup at $p\to\infty$ and using \eqref{GCap} yields 
\be
\le|\ln \frac {{\rm Cap}_{z_\infty}(\wt \K)}{{\rm Cap}_{z_\infty}( \K)} \ri|\leq  C \sqrt{\epsilon},
\ee
which implies the Hausdorff continuity of ${\rm Cap}_{z_\infty}$ over $\mathbb K(\delta, \rho)$ as required. Note that the local coordinate $z_\infty$ does not play any role since the ratio of capacities is independent of it, see \eqref{capzw}.

 \end{document}